
\catcode'32=9
\magnification=1200

\voffset=1cm

\font\tenpc=cmcsc10

\font\eightrm=cmr8
\font\eighti=cmmi8
\font\eightsy=cmsy8
\font\eightbf=cmbx8
\font\eighttt=cmtt8
\font\eightit=cmti8
\font\eightsl=cmsl8
\font\sixrm=cmr6
\font\sixi=cmmi6
\font\sixsy=cmsy6
\font\sixbf=cmbx6

\skewchar\eighti='177 \skewchar\sixi='177
\skewchar\eightsy='60 \skewchar\sixsy='60

\font\tengoth=eufm10
\font\tenbboard=msbm10
\font\eightgoth=eufm7 at 8pt
\font\eightbboard=msbm7 at 8pt
\font\sevengoth=eufm7
\font\sevenbboard=msbm7
\font\sixgoth=eufm5 at 6 pt
\font\fivegoth=eufm5

\font\tengoth=eufm10
\font\tenbboard=msbm10
\font\eightgoth=eufm7 at 8pt
\font\eightbboard=msbm7 at 8pt
\font\sevengoth=eufm7
\font\sevenbboard=msbm7
\font\sixgoth=eufm5 at 6 pt
\font\fivegoth=eufm5

\newfam\gothfam
\newfam\bboardfam

\catcode`\@=11

\def\raggedbottom{\topskip 10pt plus 36pt
\r@ggedbottomtrue}
\def\pc#1#2|{{\bigf@ntpc #1\penalty
\@MM\hskip\z@skip\smallf@ntpc #2}}

\def\tenpoint{%
  \textfont0=\tenrm \scriptfont0=\sevenrm \scriptscriptfont0=\fiverm
  \def\rm{\fam\z@\tenrm}%
  \textfont1=\teni \scriptfont1=\seveni \scriptscriptfont1=\fivei
  \def\oldstyle{\fam\@ne\teni}%
  \textfont2=\tensy \scriptfont2=\sevensy \scriptscriptfont2=\fivesy
  \textfont\gothfam=\tengoth \scriptfont\gothfam=\sevengoth
  \scriptscriptfont\gothfam=\fivegoth
  \def\goth{\fam\gothfam\tengoth}%
  \textfont\bboardfam=\tenbboard \scriptfont\bboardfam=\sevenbboard
  \scriptscriptfont\bboardfam=\sevenbboard
  \def\bboard{\fam\bboardfam}%
  \textfont\itfam=\tenit
  \def\it{\fam\itfam\tenit}%
  \textfont\slfam=\tensl
  \def\sl{\fam\slfam\tensl}%
  \textfont\bffam=\tenbf \scriptfont\bffam=\sevenbf
  \scriptscriptfont\bffam=\fivebf
  \def\bf{\fam\bffam\tenbf}%
  \textfont\ttfam=\tentt
  \def\tt{\fam\ttfam\tentt}%
  \abovedisplayskip=12pt plus 3pt minus 9pt
  \abovedisplayshortskip=0pt plus 3pt
  \belowdisplayskip=12pt plus 3pt minus 9pt
  \belowdisplayshortskip=7pt plus 3pt minus 4pt
  \smallskipamount=3pt plus 1pt minus 1pt
  \medskipamount=6pt plus 2pt minus 2pt
  \bigskipamount=12pt plus 4pt minus 4pt
  \normalbaselineskip=12pt
  \setbox\strutbox=\hbox{\vrule height8.5pt depth3.5pt width0pt}%
  \let\bigf@ntpc=\tenrm \let\smallf@ntpc=\sevenrm
  \let\petcap=\tenpc
  \normalbaselines\rm}
\def\eightpoint{%
  \textfont0=\eightrm \scriptfont0=\sixrm \scriptscriptfont0=\fiverm
  \def\rm{\fam\z@\eightrm}%
  \textfont1=\eighti \scriptfont1=\sixi \scriptscriptfont1=\fivei
  \def\oldstyle{\fam\@ne\eighti}%
  \textfont2=\eightsy \scriptfont2=\sixsy \scriptscriptfont2=\fivesy
  \textfont\gothfam=\eightgoth \scriptfont\gothfam=\sixgoth
  \scriptscriptfont\gothfam=\fivegoth
  \def\goth{\fam\gothfam\eightgoth}%
  \textfont\bboardfam=\eightbboard \scriptfont\bboardfam=\sevenbboard
  \scriptscriptfont\bboardfam=\sevenbboard
  \def\bboard{\fam\bboardfam}%
  \textfont\itfam=\eightit
  \def\it{\fam\itfam\eightit}%
  \textfont\slfam=\eightsl
  \def\sl{\fam\slfam\eightsl}%
  \textfont\bffam=\eightbf \scriptfont\bffam=\sixbf
  \scriptscriptfont\bffam=\fivebf
  \def\bf{\fam\bffam\eightbf}%
  \textfont\ttfam=\eighttt
  \def\tt{\fam\ttfam\eighttt}%
  \abovedisplayskip=9pt plus 2pt minus 6pt
  \abovedisplayshortskip=0pt plus 2pt
  \belowdisplayskip=9pt plus 2pt minus 6pt
  \belowdisplayshortskip=5pt plus 2pt minus 3pt
  \smallskipamount=2pt plus 1pt minus 1pt
  \medskipamount=4pt plus 2pt minus 1pt
  \bigskipamount=9pt plus 3pt minus 3pt
  \normalbaselineskip=9pt
  \setbox\strutbox=\hbox{\vrule height7pt depth2pt width0pt}%
  \let\bigf@ntpc=\eightrm \let\smallf@ntpc=\sixrm
  \normalbaselines\rm}

\tenpoint

\frenchspacing


\newif\ifpagetitre
\newtoks\auteurcourant \auteurcourant={\hfil}
\newtoks\titrecourant \titrecourant={\hfil}

\def\appeln@te{}
\def\vfootnote#1{\def\@parameter{#1}\insert\footins\bgroup\eightpoint
  \interlinepenalty\interfootnotelinepenalty
  \splittopskip\ht\strutbox 
  \splitmaxdepth\dp\strutbox \floatingpenalty\@MM
  \leftskip\z@skip \rightskip\z@skip
  \ifx\appeln@te\@parameter\indent \else{\noindent #1\ }\fi
  \footstrut\futurelet\next\fo@t}

\pretolerance=500 \tolerance=1000 \brokenpenalty=5000
\newdimen\hmargehaute \hmargehaute=0cm
\newdimen\lpage \lpage=13.3cm
\newdimen\hpage \hpage=20cm
\newdimen\lmargeext \lmargeext=1cm
\hsize=11.25cm
\vsize=18cm
\parskip 0pt
\parindent=12pt

\def\margehaute{\vbox to \hmargehaute{\vss}}%
\def\margebasse{\vss}

\output{\shipout\vbox to \hpage{\margehaute\nointerlineskip
  \corpsdepage\margebasse}
  \advancepageno \global\pagetitrefalse
  \ifnum\outputpenalty>-20000 \else\dosupereject\fi}

\def\corpsdepage{\hbox to \lpage{\hss\pagetexte\hskip\lmargeext}}
\def\pagetexte{\vbox{\makeheadline\pagebody\makefootline}}
\headline={\ifpagetitre\titleheadline \else
  \ifodd\pageno\rightheadline \else\leftheadline\fi\fi}
\def\leftheadline{\eightpoint\hfil\the\auteurcourant\hfil}
\def\rightheadline{\eightpoint\hfil\the\titrecourant\hfil}
\def\titleheadline{\hfill}
\pagetitretrue

\def\footnoterule{\kern-6\p@
  \hrule width 2truein \kern 5.6\p@} 

\def\pd#1#2 {\pc#1#2| }

\def\pointir{\discretionary{.}{}{.\kern.35em---\kern.7em}\nobreak
\hskip 0em plus .3em minus .4em }

\def\abstract#1{\vbox{\eightpoint \pc ABSTRACT|\pointir #1}}

\def\titre#1|{\message{#1}
              \par\vskip 30pt plus 24pt minus 3pt\penalty -1000
              \vskip 0pt plus -24pt minus 3pt\penalty -1000
              \centerline{\bf #1}
              \vskip 5pt
              \penalty 10000 }

\def\section#1|{\par\vskip .3cm
                {\bf #1}.\quad}

\def\ssection#1|{\par\vskip .2cm
                {\it #1}\pointir}

\long\def\th#1|#2\finth{\par\medskip
              {\petcap #1\pointir}{\it #2}\par\smallskip}

\long\def\tha#1|#2\fintha{\par\medskip
                    {\petcap #1.}\par\nobreak{\it #2}\par\smallskip}

\def\rem#1|{\par\medskip
            {{\it #1}\pointir}}

\def\rema#1|{\par\medskip
             {{\it #1.}\par\nobreak }}

\def\article#1|#2|#3|#4|#5|#6|#7|
    {{\leftskip=7mm\noindent
     \hangindent=2mm\hangafter=1
     \llap{[#1]\hskip.35em}{#2}.\quad
     #3, {\sl #4}, vol.\nobreak\ {\bf #5}, {\oldstyle #6},
     p.\nobreak\ #7.\par}}
\def\livre#1|#2|#3|#4|
    {{\leftskip=7mm\noindent
    \hangindent=2mm\hangafter=1
    \llap{[#1]\hskip.35em}{#2}.\quad
    {\sl #3}.\quad #4.\par}}
\def\divers#1|#2|#3|
    {{\leftskip=7mm\noindent
    \hangindent=2mm\hangafter=1
     \llap{[#1]\hskip.35em}{#2}.\quad
     #3.\par}}
\mathchardef\conj="0365
\def\proof{\par{\it Proof}.\quad}
\def\qed{\quad\raise -2pt\hbox{\vrule\vbox to 10pt{\hrule width 4pt
\vfill\hrule}\vrule}}

\def\cqfd{\penalty 500 \hbox{\qed}\par\smallskip}
\def\decale#1|{\par\noindent\hskip 28pt\llap{#1}\kern 5pt}

\catcode`\@=12

\def\Grille{\setbox13=\vbox to 5mm{\hrule width 110mm\vfill}
\setbox13=\vbox{\offinterlineskip
   \copy13\copy13\copy13\copy13\copy13\copy13\copy13\copy13
   \copy13\copy13\copy13\copy13\box13\hrule width 110mm}
\setbox14=\hbox to 5mm{\vrule height 65mm\hfill}
\setbox14=\hbox{\copy14\copy14\copy14\copy14\copy14\copy14
   \copy14\copy14\copy14\copy14\copy14\copy14\copy14\copy14
   \copy14\copy14\copy14\copy14\copy14\copy14\copy14\copy14\box14}
\ht14=0pt\dp14=0pt\wd14=0pt
\setbox13=\vbox to 0pt{\vss\box13\offinterlineskip\box14}
\wd13=0pt\box13}


\def\fleche(#1,#2)\dir(#3,#4)\long#5{%
\noalign{\nointerlineskip\leftput(#1,#2){\vector(#3,#4){#5}}\nointerlineskip}}


\def\hfl#1#2#3{\smash{\mathop{\hbox to#3{\rightarrowfill}}\limits
^{\scriptstyle#1}_{\scriptstyle#2}}}

\def\gfl#1#2#3{\smash{\mathop{\hbox to#3{\leftarrowfill}}\limits
^{\scriptstyle#1}_{\scriptstyle#2}}}


 \message{`lline' & `vector' macros from LaTeX}
 \catcode`@=11
\def\{{\relax\ifmmode\lbrace\else$\lbrace$\fi}
\def\}{\relax\ifmmode\rbrace\else$\rbrace$\fi}
\def\newcount{\alloc@0\count\countdef\insc@unt}
\def\newdimen{\alloc@1\dimen\dimendef\insc@unt}
\def\newwrite{\alloc@7\write\chardef\sixt@@n}

\newwrite\@unused
\newcount\@tempcnta
\newcount\@tempcntb
\newdimen\@tempdima
\newdimen\@tempdimb
\newbox\@tempboxa

\def\@spaces{\space\space\space\space}
\def\@whilenoop#1{}
\def\@whiledim#1\do #2{\ifdim #1\relax#2\@iwhiledim{#1\relax#2}\fi}
\def\@iwhiledim#1{\ifdim #1\let\@nextwhile=\@iwhiledim
        \else\let\@nextwhile=\@whilenoop\fi\@nextwhile{#1}}
\def\@badlinearg{\@latexerr{Bad \string\line\space or \string\vector
   \space argument}}
\def\@latexerr#1#2{\begingroup
\edef\@tempc{#2}\expandafter\errhelp\expandafter{\@tempc}%
\def\@eha{Your command was ignored.
^^JType \space I <command> <return> \space to replace it
  with another command,^^Jor \space <return> \space to continue without it.}
\def\@ehb{You've lost some text. \space \@ehc}
\def\@ehc{Try typing \space <return>
  \space to proceed.^^JIf that doesn't work, type \space X <return> \space to
  quit.}
\def\@ehd{You're in trouble here.  \space\@ehc}

\typeout{LaTeX error. \space See LaTeX manual for explanation.^^J
 \space\@spaces\@spaces\@spaces Type \space H <return> \space for
 immediate help.}\errmessage{#1}\endgroup}
\def\typeout#1{{\let\protect\string\immediate\write\@unused{#1}}}

\font\tenln    = line10
\font\tenlnw   = linew10

\newdimen\@wholewidth
\newdimen\@halfwidth
\newdimen\unitlength 

\unitlength =1pt


\def\thinlines{\let\@linefnt\tenln \let\@circlefnt\tencirc
  \@wholewidth\fontdimen8\tenln \@halfwidth .5\@wholewidth}
\def\thicklines{\let\@linefnt\tenlnw \let\@circlefnt\tencircw
  \@wholewidth\fontdimen8\tenlnw \@halfwidth .5\@wholewidth}

\def\linethickness#1{\@wholewidth #1\relax \@halfwidth .5\@wholewidth}

\newif\if@negarg

\def\lline(#1,#2)#3{\@xarg #1\relax \@yarg #2\relax
\@linelen=#3\unitlength
\ifnum\@xarg =0 \@vline
  \else \ifnum\@yarg =0 \@hline \else \@sline\fi
\fi}

\def\@sline{\ifnum\@xarg< 0 \@negargtrue \@xarg -\@xarg \@yyarg -\@yarg
  \else \@negargfalse \@yyarg \@yarg \fi
\ifnum \@yyarg >0 \@tempcnta\@yyarg \else \@tempcnta -\@yyarg \fi
\ifnum\@tempcnta>6 \@badlinearg\@tempcnta0 \fi
\setbox\@linechar\hbox{\@linefnt\@getlinechar(\@xarg,\@yyarg)}%
\ifnum \@yarg >0 \let\@upordown\raise \@clnht\z@
   \else\let\@upordown\lower \@clnht \ht\@linechar\fi
\@clnwd=\wd\@linechar
\if@negarg \hskip -\wd\@linechar \def\@tempa{\hskip -2\wd\@linechar}\else
     \let\@tempa\relax \fi
\@whiledim \@clnwd <\@linelen \do
  {\@upordown\@clnht\copy\@linechar
   \@tempa
   \advance\@clnht \ht\@linechar
   \advance\@clnwd \wd\@linechar}%
\advance\@clnht -\ht\@linechar
\advance\@clnwd -\wd\@linechar
\@tempdima\@linelen\advance\@tempdima -\@clnwd
\@tempdimb\@tempdima\advance\@tempdimb -\wd\@linechar
\if@negarg \hskip -\@tempdimb \else \hskip \@tempdimb \fi
\multiply\@tempdima \@m
\@tempcnta \@tempdima \@tempdima \wd\@linechar \divide\@tempcnta \@tempdima
\@tempdima \ht\@linechar \multiply\@tempdima \@tempcnta
\divide\@tempdima \@m
\advance\@clnht \@tempdima
\ifdim \@linelen <\wd\@linechar
   \hskip \wd\@linechar
  \else\@upordown\@clnht\copy\@linechar\fi}

\def\@hline{\ifnum \@xarg <0 \hskip -\@linelen \fi
\vrule height \@halfwidth depth \@halfwidth width \@linelen
\ifnum \@xarg <0 \hskip -\@linelen \fi}

\def\@getlinechar(#1,#2){\@tempcnta#1\relax\multiply\@tempcnta 8
\advance\@tempcnta -9 \ifnum #2>0 \advance\@tempcnta #2\relax\else
\advance\@tempcnta -#2\relax\advance\@tempcnta 64 \fi
\char\@tempcnta}

\def\vector(#1,#2)#3{\@xarg #1\relax \@yarg #2\relax
\@linelen=#3\unitlength
\ifnum\@xarg =0 \@vvector
  \else \ifnum\@yarg =0 \@hvector \else \@svector\fi
\fi}

\def\@hvector{\@hline\hbox to 0pt{\@linefnt
\ifnum \@xarg <0 \@getlarrow(1,0)\hss\else
    \hss\@getrarrow(1,0)\fi}}

\def\@vvector{\ifnum \@yarg <0 \@downvector \else \@upvector \fi}

\def\@svector{\@sline
\@tempcnta\@yarg \ifnum\@tempcnta <0 \@tempcnta=-\@tempcnta\fi
\ifnum\@tempcnta <5
  \hskip -\wd\@linechar
  \@upordown\@clnht \hbox{\@linefnt  \if@negarg
  \@getlarrow(\@xarg,\@yyarg) \else \@getrarrow(\@xarg,\@yyarg) \fi}%
\else\@badlinearg\fi}

\def\@getlarrow(#1,#2){\ifnum #2 =\z@ \@tempcnta='33\else
\@tempcnta=#1\relax\multiply\@tempcnta \sixt@@n \advance\@tempcnta
-9 \@tempcntb=#2\relax\multiply\@tempcntb \tw@
\ifnum \@tempcntb >0 \advance\@tempcnta \@tempcntb\relax
\else\advance\@tempcnta -\@tempcntb\advance\@tempcnta 64
\fi\fi\char\@tempcnta}

\def\@getrarrow(#1,#2){\@tempcntb=#2\relax
\ifnum\@tempcntb < 0 \@tempcntb=-\@tempcntb\relax\fi
\ifcase \@tempcntb\relax \@tempcnta='55 \or
\ifnum #1<3 \@tempcnta=#1\relax\multiply\@tempcnta
24 \advance\@tempcnta -6 \else \ifnum #1=3 \@tempcnta=49
\else\@tempcnta=58 \fi\fi\or
\ifnum #1<3 \@tempcnta=#1\relax\multiply\@tempcnta
24 \advance\@tempcnta -3 \else \@tempcnta=51\fi\or
\@tempcnta=#1\relax\multiply\@tempcnta
\sixt@@n \advance\@tempcnta -\tw@ \else
\@tempcnta=#1\relax\multiply\@tempcnta
\sixt@@n \advance\@tempcnta 7 \fi\ifnum #2<0 \advance\@tempcnta 64 \fi
\char\@tempcnta}

\def\@vline{\ifnum \@yarg <0 \@downline \else \@upline\fi}

\def\@upline{\hbox to \z@{\hskip -\@halfwidth \vrule
  width \@wholewidth height \@linelen depth \z@\hss}}

\def\@downline{\hbox to \z@{\hskip -\@halfwidth \vrule
  width \@wholewidth height \z@ depth \@linelen \hss}}

\def\@upvector{\@upline\setbox\@tempboxa\hbox{\@linefnt\char'66}\raise
     \@linelen \hbox to\z@{\lower \ht\@tempboxa\box\@tempboxa\hss}}

\def\@downvector{\@downline\lower \@linelen
      \hbox to \z@{\@linefnt\char'77\hss}}

\thinlines

\newcount\@xarg
\newcount\@yarg
\newcount\@yyarg
\newcount\@multicnt
\newdimen\@xdim
\newdimen\@ydim
\newbox\@linechar
\newdimen\@linelen
\newdimen\@clnwd
\newdimen\@clnht
\newdimen\@dashdim
\newbox\@dashbox
\newcount\@dashcnt
 \catcode`@=12


\newbox\tbox
\newbox\tboxa

\def\leftzer#1{\setbox\tbox=\hbox to 0pt{#1\hss}%
     \ht\tbox=0pt \dp\tbox=0pt \box\tbox}

\def\rightzer#1{\setbox\tbox=\hbox to 0pt{\hss #1}%
     \ht\tbox=0pt \dp\tbox=0pt \box\tbox}

\def\centerzer#1{\setbox\tbox=\hbox to 0pt{\hss #1\hss}%
     \ht\tbox=0pt \dp\tbox=0pt \box\tbox}

%
\def\image(#1,#2)#3{\vbox to #1{\offinterlineskip
    \vss #3 \vskip #2}}


\def\leftput(#1,#2)#3{\setbox\tboxa=\hbox{%
    \kern #1\unitlength
    \raise #2\unitlength\hbox{\leftzer{#3}}}%
    \ht\tboxa=0pt \wd\tboxa=0pt \dp\tboxa=0pt\box\tboxa}

\def\rightput(#1,#2)#3{\setbox\tboxa=\hbox{%
    \kern #1\unitlength
    \raise #2\unitlength\hbox{\rightzer{#3}}}%
    \ht\tboxa=0pt \wd\tboxa=0pt \dp\tboxa=0pt\box\tboxa}

\def\centerput(#1,#2)#3{\setbox\tboxa=\hbox{%
    \kern #1\unitlength
    \raise #2\unitlength\hbox{\centerzer{#3}}}%
    \ht\tboxa=0pt \wd\tboxa=0pt \dp\tboxa=0pt\box\tboxa}

\unitlength=1mm

\def\put(#1,#2)#3{\noalign{\nointerlineskip
                               \centerput(#1,#2){$#3$}
                                \nointerlineskip}}
\def\fleche(#1,#2)\dir(#3,#4)\long#5{%
{\leftput(#1,#2){\vector(#3,#4){#5}}}}


\catcode`\@=11
\def\matrice#1{\null \,\vcenter {\normalbaselines \m@th
\ialign {\hfil $##$\hfil &&\  \hfil $##$\hfil\crcr
\mathstrut \crcr \noalign {\kern -\baselineskip } #1\crcr
\mathstrut \crcr \noalign {\kern -\baselineskip }}}\,}

\def\petitematrice#1{\left(\null\vcenter {\normalbaselines \m@th
\ialign {\hfil $##$\hfil &&\thinspace  \hfil $##$\hfil\crcr
\mathstrut \crcr \noalign {\kern -\baselineskip } #1\crcr
\mathstrut \crcr \noalign {\kern -\baselineskip }}}\right)}

\catcode`\@=12


\def\negg{\mathop{\rm neg}\nolimits}
\def\inv{\mathop{\rm inv}\nolimits}
\def\imaj{\mathop{\rm imaj}\nolimits}
\def\tot{\mathop{\rm tot}\nolimits}
\def\L{\mathop{\rm L\kern 0pt}}

\def\odd{\mathop{\rm odd}\nolimits}

\def\NIW{\mathop{{\hbox{\eightrm NIW}}}\nolimits}
\def\DW{\mathop{{\hbox{\eightrm DW}}}\nolimits}

\def\WSP{\mathop{{\hbox{\eightrm WSP}}}\nolimits}
\def\WSD{\mathop{{\hbox{\eightrm WSD}}}\nolimits}
\def\fix{\mathop{\rm fix}\nolimits}
\def\pix{\mathop{\rm pix}\nolimits}

\def\fmaj{\mathop{\rm fmaj}\nolimits}
\def\fdes{\mathop{\rm fdes}\nolimits}
\def\des{\mathop{\rm des}\nolimits}
\def\maj{\mathop{\rm maj}\nolimits}

\def\Ligne{\mathop{\rm Ligne}\nolimits}
\def\Iligne{\mathop{\rm Iligne}\nolimits}

\def\odd{\mathop{\rm odd}\nolimits}

\def\Fix{\mathop{\rm Fix}\nolimits}
\def\Pix{\mathop{\rm Pix}\nolimits}
\def\negg{\mathop{\rm neg}\nolimits}
\def\Neg{\mathop{\rm Neg}\nolimits}

\titrecourant={FIXED AND PIXED POINTS}
\auteurcourant={DOMINIQUE FOATA AND GUO-NIU HAN}

\rightline{2006/06/23}
\vglue .5cm

\centerline{\bf SIGNED WORDS AND PERMUTATIONS, IV;}
\smallskip
\centerline{\bf FIXED AND PIXED POINTS}
\bigskip
\centerline{
\bf Dominique Foata and Guo-Niu Han}

\bigskip
\hbox{\hskip5cm\vbox{\eightpoint
\eightsl\vbox{\halign{#\hfill\cr
Von Jacobs hat er die Statur,\cr
Des Rechnens ernstes F\"uhren,\cr
Von Lott\"archen die Frohnatur\cr
und Lust zu diskretieren.\cr
\noalign{\medskip}
To Volker Strehl, a dedication \`a la Goethe,\cr
on the occasion of his sixtieth birthday.\cr}}}}

\bigskip
\centerline{\bf Abstract}
\smallskip
{\narrower\eightpoint

\noindent
The flag-major index ``fmaj" and the classical length function
``$\ell$" are used to construct two $q$-analogs of the
generating polynomial for the hyperoctahedral group~$B_n$
by number of positive and negative fixed points (resp. pixed
points). Specializations of those $q$-analogs are also
derived dealing with signed derangements and desarrangements,
as well as several classical results that were previously proved
for the symmetric group.

\bigskip}

\centerline{\bf 1. Introduction}

\medskip
The\footnote{}{2000 {\it Mathematics Subject Classification.}
Primary 05A15, 05A30, 05E15.\hfil\break\indent
{\it Key words and phrases.} Hyperoctahedral group, length
function, flag-major index, signed permutations, fixed points,
pixed points, derangements, desarrangements, pixed factorization.}
statistical study of the hyperoctahedral group~$B_n$, initiated by
Reiner ([Re93a], [Re93b], [Re93c], [Re95a], [Re95b]), has been rejuvenated
by Adin and Roichman [AR01] with their introduction of the
{\it flag-major index}, which was shown [ABR01] to be equidistributed
with the {\it length function}. See also their recent papers on the
subject [ABR05], [ReRo05]. It then appeared natural to extend the
numerous results obtained for the symmetric group~${\goth S}_n$ to
the groug~$B_n$. Furthermore, flag-major index and length
function become the true $q$-analog makers needed for
calculating various multivariable distributions on~$B_n$.

In the present paper we start with a generating polynomial
for~$B_n$ by a three-variable statistic involving the number of
fixed points (see formula (1.3)) and show that there are two ways of
$q$-analogizing it, by using the flag-major index on the one
hand, and the length function, on the other hand. As will be
indicated, the introduction of an extra variable~$Z$ makes it
possible to specialize all our results to the symmetric group. Let
us first give the necessary notations.

Let $B_n$ be the hyperoctahedral group of all {\it signed
permutations} of order~$n$. The elements of~$B_n$ may be viewed
as words $w=x_1x_2\cdots x_n$, where each~$x_i$ belongs to
$\{-n,\ldots, -1,1,\ldots,n\}$ and $|x_1||x_2|\cdots |x_n|$ is a
permutation of $12\ldots n$. The {\it set} (resp. the {\it number}) of
{\it negative} letters among the~$x_i$'s is denoted by $\Neg w$ (resp.
$\negg w$).  A {\it positive fixed point} of the signed permutation
$w=x_1x_2\cdots x_n$ is a (positive) integer~$i$ such that
$x_i=i$. It is convenient to write $\overline i:=-i$ for each
integer~$i$. Also, when~$A$ is a set of integers, let
$\overline A:=\{\overline i:i\in A\}$. If $x_i=\overline i$ with $i$
positive, we say that~$\overline i$ is a {\it negative fixed point}
of~$w$. The set of all positive (resp. negative) fixed points of~$w$ is
denoted by $\Fix^+w$ (resp. $\Fix^- w$). Notice that
$\Fix^-w\subset \Neg w$. Also let
$$
\fix^+w:=\#\Fix^+w;\quad \fix^-w:=\#\Fix^-w.\leqno(1.1)
$$
There are $2^nn!$ signed permutations of order~$n$.
The symmetric group~${\goth S}_n$ may be considered as the
subset of all~$w$ from~$B_n$ such that $\Neg w=\emptyset$.

\medskip
The purpose of this paper is to provide {\it two $q$-analogs}
for the polynomials $B_n(Y_0,Y_1,Z)$ defined by the identity
$$
\sum_{n\ge 0}{u^n\over n!}B_n(Y_0,Y_1,Z)
=\bigl(1-u(1+Z)\bigr)^{-1}
\times {\exp(u(Y_0+Y_1Z))\over \exp(u(1+Z))}.\leqno(1.2)
$$
When $Z=0$, the right-hand side becomes
$(1-u)^{-1}\exp(uY_0)/\exp(u)$, which is the exponential
generating function for the generating polynomials for the
groups~${\goth S}_n$ by number of fixed points (see [Ri58],
chap.~4). Also, by identification, 
$B_n(1,1,1)=2^nn!$ and it is easy to show (see Theorem~1.1) that
$B_n(Y_0,Y_1,Z)$ is in fact the generating polynomial for the
group~$B_n$ by the three-variable statistic
$(\fix^+,\fix^-,\negg)$, that is,
$$
B_n(Y_0,Y_1,Z)=\sum_{w\in
B_n}Y_0^{\fix^+\!w}\,Y_1^{\fix^-\!w}\, Z^{\negg
w}.\leqno(1.3)
$$

Recall the traditional notations for the $q$-ascending
factorials
$$\displaylines{
\rlap{(1.4)}\hfill
(a;q)_n:=\cases{1,&if $n=0$;\cr
(1-a)(1-aq)\cdots (1-aq^{n-1}),&if $n\ge 1$;\cr}\hfill\cr
(a;q)_\infty:=\prod_{n\ge 1}(1-aq^{n-1});\kern4.35cm\cr
\noalign{\hbox{for the $q$-multinomial coefficients}}
\rlap{(1.5)}\hfill
{\,n\,\brack m_1,\ldots,m_k}_q:={(q;q)_n\over
(q;q)_{m_1}\cdots (q;q)_{m_k}}
\quad (m_1+\cdots+m_k= n);\hfill\cr
\noalign{\hbox{and for the two $q$-exponentials
(see [GaRa90, chap.~1])}}
\rlap{(1.6)}\hfill
e_q(u)=\sum_{n\ge 0}{u^n\over (q;q)_n}={1\over
(u;q)_\infty};\quad 
E_q(u)=\sum_{n\ge 0}{q^{n\choose2}u^n\over (q;q)_n}=(-u;q)_\infty.\cr}
$$

\goodbreak

Our two $q$-analogs, denoted by
${}^{\ell\kern-1.5pt}B_n(q,Y_0,Y_1,Z)$ and $B_n(q,Y_0,Y_1,Z)$, are
respectively defined by the identities:
$$\displaylines{{(1.7)}\quad
\sum_{n\ge 0}
{u^n\over (-Zq;q)_n\,(q;q)_n}\,{}^{\ell\kern-1.5pt}B_n(q,Y_0,Y_1,Z)\hfill\cr
\noalign{\vskip-5pt}
\hfill{}=\Bigl(1-{u\over 1-q}\Bigr)^{-1}\times (u;q)_\infty\,
\Bigl(\sum_{n\ge 0}
{(-qY_0^{-1}Y_1Z;q)_n\,(uY_0)^n\over
(-Zq;q)_n\,(q;q)_n}\Bigr);
\quad\cr
(1.8)\quad
\sum_{n\ge 0}{u^n\over
(q^2;q^2)_n}B_n(q,Y_0,Y_1,Z)\hfill\cr
\noalign{\vskip-10pt}
\hfill{} =\Bigl(1-u{1+qZ\over
1-q^2}\Bigr)^{-1}
\kern-4pt\times\kern-2pt
{(u;q^2)_\infty  \over (uY_0;q^2)_\infty}
{(-uqY_1Z;q^2)_\infty\over (-uqZ;q^2)_\infty}
.\kern15pt\cr}
$$
Those two identities can be shown to yield (1.2) when
$q=1$.

There is also a graded form of (1.8) in the sense that an
extra variable~$t$ can be added to form a new polynomial
$B_n(t,q,Y_0,Y_1,Z)$ with nonnegative integral coefficients
that specializes into $B_n(q,Y_0,Y_1,Z)$ for $t=1$. Those
polynomials are defined by the identity
$$
\displaylines{(1.9)\quad
\sum_{n\ge 0}(1+t)B_n(t,q,Y_0,Y_1,Z){u^n\over
(t^2;q^2)_{n+1}}
\hfill\cr
\noalign{\vskip-5pt}
\hfill{}
=\sum_{s\ge 0} t^s
\Bigl(1-u\sum_{i=0}^s q^iZ\strut^{\chi (i\ {\rm
odd})}\Bigr)^{-1}\kern-4pt\times\kern-2pt
{(u;q^2)_{\lfloor s/2\rfloor+1}  \over 
(uY_0;q^2)_{\lfloor s/2\rfloor+1}}
{(-uqY_1Z;q^2)_{\lfloor (s+1)/2\rfloor}\over
(-uqZ;q^2)_{\lfloor (s+1)/2\rfloor}},\cr}
$$
where for each statement~$A$ we let $\chi(A)=1$ or~0
depending on whether~$A$ is true or not. The importance of identity (1.9) lies in its
numerous specializations, as can be seen in Fig.~1.  

\medskip
The two $q$-extensions ${}^{\ell\kern-1.5pt}B_n(q,Y_0,Y_1,Z)$ and
$B_n(t,q,Y_0,Y_1,Z)$ being now defined, the program is to
derive appropriate combinatorial interpretations for them.
Before doing so we need have a second combinatorial interpretation
for the polynomial $B_n(Y_0,Y_1,Z)$ besides the one mentioned in
(1.3). Let $w=x_1x_2\cdots x_n$ be a word, all letters of which are
integers without any repetitions. Say that~$w$ is a {\it
desarrangement} if $x_1>x_2>\cdots >x_{2k}$ and
$x_{2k}<x_{2k+1}$ for some~$k\ge 1$. By convention,
$x_{n+1}=\infty$. We could also say that the {\it leftmost
trough} of~$w$ occurs at an {\it even} position. This notion
was introduced by D\'esarm\'enien [De84] and elegantly used
in a subsequent paper [DeWa88]. Notice that there is no one-letter
desarrangement. By convention, the empty word~$e$ is also a
desarrangement.

Now let $w=x_1x_2\cdots x_n$ be a signed permutation.
Unless~$w$ is increasing, there is always a nonempty right
factor of~$w$ which is a desarrangement. It then makes sense to
define~$w^d$ as the {\it longest} such a right factor.
Hence,~$w$ admits a unique factorization $w=w^-w^+w^d$,
called its {\it pixed}\footnote{$^{(1)}$}{``Pix," of course, must not
be taken here for the abbreviated form of ``pictures."}
{\it factorization}, where $w^-$ and $w^+$ are both  {\it
increasing}, the letters of~$w^-$ being {\it negative}, those of
$w^+$ {\it positive} and where~$w^d$ is the longest right factor
of~$w$ which is a desarrangement.

\medskip
For example, the pixed factorizations of the following signed
permutations are materialized by vertical bars: 
$w=\overline 5\,\overline2\mid e\mid \overline 3\,\overline
4\,1$; $w=\overline
5\mid e\mid\overline2\,\overline 3\,1\,\overline 4$;\quad
$w=\overline 5\,\overline3\,\overline 2\mid 1\, 4\mid e$;\quad 
$w=\overline 5\,\overline3\mid 1\mid 4\,2$;\quad
$w=\overline 5\,\overline3\mid e\mid 4\, 1\,2$.

\medskip
Let $w=w^-w^+w^d$ be the pixed factorization of 
$w=x_1x_2\cdots x_n$.
If $w^-=x_1\cdots x_k$,
$w^+=x_{k+1}\cdots x_{k+l}$, define
$\Pix^-w:=\{x_1,\ldots,x_k\}$,
$\Pix^+w:=\{x_{k+1},\ldots, x_{k+l}\}$,
$\pix^-w:=\#\Pix^-w$, $\pix^+w:=\#\Pix^+w$.

\proclaim Theorem 1.1. The polynomial 
$B_n(Y_0,Y_1,Z)$ defined by $(1.2)$ admits the following two
combinatorial interpretations:
$$B_n(Y_0,Y_1,Z)
=\!\!\sum_{w\in B_n}Y_1^{\fix^-w}\,Y_2^{\fix^+w}
\,Z^{\negg w}
=\!\!\sum_{w\in B_n}Y_1^{\pix^-w}\,Y_2^{\pix^+w}
\,Z^{\negg w}.
$$

Theorem~1.1 is proved in section~2. A 
bijection~$\phi$ of~$B_n$ onto itself will be constructed that
satisfies 
$(\Fix^-,\Fix^+,\Neg)
\,w=(\Pix^-,\Pix^+,\Neg)\,\phi(w)$.

\medskip
Let ``$\ell$" be the length function of~$B_n$ (see [Bo68, p.~7],
[Hu90, p.~12] or the working definition given in (3.1)). 
As seen in Theorem~1.2, 
``$\ell$" is to be added to the three-variable statistic
$(\pix^+,\pix^-,\negg)$ (and not to
$(\fix^+,\fix^-,\negg)$) for deriving the combinatorial
interpretation of
${}^{\ell\kern-1.5pt}B_n(q,Y_0,Y_1,Z)$. This theorem is proved in
section~3.

\proclaim Theorem 1.2. For each $n\ge 0$ let 
${}^{\ell\kern-1.5pt}B_n(q,Y_0,Y_1,Z)$ be the polynomial defined in
$(1.7)$. Then
$$
{}^{\ell\kern-1.5pt}B_n(q,Y_0,Y_1,Z)=\sum_{w\in B_n}q^{\ell(w)}\,
Y_0^{\pix^+w}\,Y_1^{\pix^- w}\,
Z^{\negg
w}.\leqno(1.10)
$$

The variables $t$ and $q$ which are added to
interpret our second extension $B_n(t,q,Y_0,Y_1,Z)$ will carry the
flag-descent number ``fdes" and the flag-major index ``fmaj."
For each signed permutation $w=x_1x_2\cdots x_n$ the usual
{\it number of descents} ``des" is defined by $\des
w:=\sum_{i=1}^{n-1}\chi(x_i>x_{i+1})$, the {\it major
index} ``maj" by $\maj w:=\sum_{i=1}^{n-1}i
\,\chi(x_i>x_{i+1})$, the {\it flag descent number} ``$\fdes $"
and the {\it flag-major index} ``$\fmaj $" by
$$
\fdes w:=2\des w+\chi(x_1<0);\quad
\fmaj w:=2\maj w+\negg w.\leqno(1.11)
$$

\proclaim Theorem 1.3. For each $n\ge 0$ let 
$B_n(t,q,Y_0,Y_1,Z)$ be the polynomial defined in
$(1.9)$. Then
$$
B_n(t,q,Y_0,Y_1,Z)=\sum_{w\in B_n}t^{\fdes w}\,q^{\fmaj
w}\,Y_0^{\fix^+w}\,Y_1^{\fix^- w}\,Z^{\negg w}.
\leqno(1.12)
$$

Theorem 1.3 is proved in Section~5 after discussing the
combinatorics of the so-called {\it weighted signed
permutations} in Section~4. Section~6 deals with
numerous specializations of Theorem~1.2 and~1.3 obtained by
taking numerical values, essentially~0 or~1, for certain
variables. Those specializations are illustrated by the
following diagram (Fig.~1). When $Z=0$, the statistic ``neg"
plays no role and the signed permutations become plain
permutations; the second column of the diagram is then
mapped on the third one that only involves generating
polynomials for~${\goth S}_n$ or subsets of that group.

\midinsert
{\eightpoint

\newbox\boxdiag
\setbox\boxdiag=\vbox{\offinterlineskip
\fleche(-29,83)\dir(1,0)\long{17}
\centerput(-20,84.5){$\scriptstyle (\fmaj,\fix^-,\negg)$}
\centerput(-33,83){$D_n^B$}
\centerput(0,83){$D_n^B(q,Y_1,Z)$}
\fleche(0,82)\dir(0,-1)\long{7}
\centerput(3,78){$\scriptstyle \fdes$}
\centerput(30,83){$D_n(q)$}
\fleche(30,82)\dir(0,-1)\long{7}
\centerput(33,78){$\scriptstyle \des$}
\centerput(18,84.5){$\scriptstyle \negg$}
\fleche(21,83)\dir(-1,0)\long{7}
\fleche(43.5,83)\dir(-1,0)\long{7}
\centerput(40,84.5){$\scriptstyle\maj$}
\centerput(48,83){$D_n$}
\centerput(0,72){$D_n^B(t,q,Y_1,Z)$}
\fleche(0,70)\dir(0,-1)\long{7}
\centerput(3.5,66){$\scriptstyle \fix^+$}
\centerput(18,73.5){$\scriptstyle \negg$}
\fleche(21,72)\dir(-1,0)\long{7}
\fleche(30,70)\dir(0,-1)\long{7}
\centerput(30,72){$D_n(t,q)$}
\centerput(33.5,66){$\scriptstyle \fix^+$}
\centerput(0,60){$\bf B_n(t,q,Y_0,Y_1,Z)$}
\centerput(18,62.5){$\scriptstyle \negg$}
\fleche(0,52)\dir(0,1)\long{7}
\centerput(3,54){$\scriptstyle\fdes$}
\centerput(0,49){$B_n(q,Y_0,Y_1,Z)$}
\centerput(30,49){$A_n(q,Y_0)$}
\fleche(30,52)\dir(0,1)\long{7}
\centerput(33,54){$\scriptstyle\des$}
\centerput(30,60){$A_n(t,q,Y_0)$}
\fleche(-29,37)\dir(1,0)\long{17}
\centerput(-20,39){$\scriptstyle (\fix^-,\fix^+,\negg)$}
\centerput(-20,34){$\scriptstyle (\pix^-,\pix^+,\negg)$}
\centerput(-32,36.5){$B_n$}
\centerput(0,36.4){$\bf B_n(Y_0,Y_1,Z)$}
\fleche(0,39.5)\dir(0,1)\long{7}
\centerput(3,42){$\scriptstyle\fmaj$}
\centerput(30,36.5){$A_n(Y_0)$}
\fleche(30,39.5)\dir(0,1)\long{7}
\centerput(33,42){$\scriptstyle\maj$}
\fleche(43.5,36.5)\dir(-1,0)\long{7}
\centerput(40,38){$\scriptstyle\fix^+$}
\centerput(40,34){$\scriptstyle\pix^+$}
\centerput(48,36){${\goth S}_n$}
\fleche(30,35)\dir(0,-1)\long{7}
\centerput(33,31){$\scriptstyle\inv$}
\centerput(0,24.5){$\bf {}^{\ell}B_n(q,Y_0,Y_1,Z)$}
\fleche(0,34.5)\dir(0,-1)\long{7}
\centerput(30,24.5){${}^{\inv\kern-1pt}A_n(q,Y_0 )$}
\centerput(1.5,31){$\scriptstyle \ell$}
\centerput(33.5,66){$\scriptstyle\fix^+$}
\fleche(21,61)\dir(-1,0)\long{7}
\fleche(0,16)\dir(0,1)\long{7}
\centerput(3,19){$\scriptstyle \pix^+$}
\centerput(0,12){${}^{\ell\kern-1.5pt}K_n^B(q,Y_1,Z)$}
\centerput(33.5,19){$\scriptstyle \pix^+$}
\centerput(30,12){${}^{\inv}K_n(q)$}
\fleche(-29,13)\dir(1,0)\long{17}
\centerput(-20,10){$\scriptstyle (\ell,\pix^-,\negg)$}
\centerput(-32,12){$K^B_n$}
\fleche(30,16)\dir(0,1)\long{7}
\fleche(43.5,13)\dir(-1,0)\long{7}
\centerput(40,14){$\scriptstyle\inv$}
\centerput(40,10){$\scriptstyle\imaj$}
\centerput(48,12){$K_n$}
\centerput(18,14.5){$\scriptstyle \negg$}
\fleche(21,13)\dir(-1,0)\long{7}
\centerput(18,26.5){$\scriptstyle \negg$}
\fleche(21,25)\dir(-1,0)\long{7}
\centerput(18,37.5){$\scriptstyle \negg$}
\fleche(21,36)\dir(-1,0)\long{7}
\centerput(18,49.5){$\scriptstyle \negg$}
\fleche(21,48)\dir(-1,0)\long{7}
}

\vglue8.6cm
\centerline{\box\boxdiag\hskip2.8cm}
\vskip-.8cm

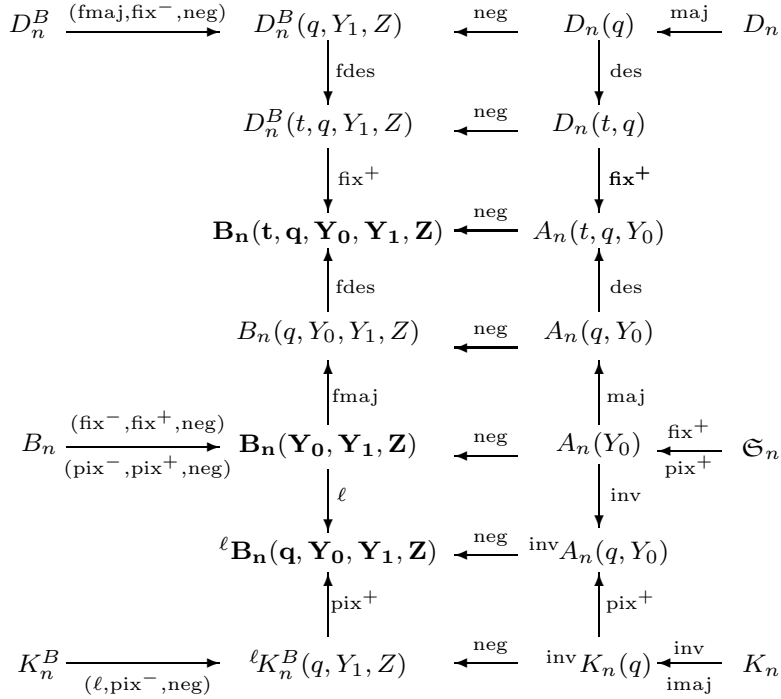
\centerline{Fig. 1}
\vskip-2pt

}

\endinsert

The first (resp. fourth) column refers to specific subsets of~$B_n$
(resp. of~${\goth S}_n$):

\vskip-24pt
$$\eqalign{
D_n&:=\{w\in B_n:\Fix^+w=\Neg w=\emptyset\};\cr
K_n&:=\{w\in B_n:\Pix^+w=\Neg w=\emptyset\};\cr
D_n^B&:=\{w\in B_n:\Fix^+w=\emptyset\};\cr
K_n^B&:=\{w\in B_n:\Pix^+w=\emptyset\}.\cr}\leqno(1.13)
$$
The elements of $D_n$ are the classical {\it derangements} and
provide the most natural combinatorial interpretations of the
derangement numbers $d_n=\#D_n$ (see [Co70], p.~9--12). By
analogy, the elements of~$D_n^B$ are called {\it signed
derangements}. They have been studied by Chow [Ch06] in a recent
note. The elements of~$K_n$ (resp. of $K_n^B$) are called {\it
desarrangements} (resp. {\it signed desarrangements}) of order~$n$.
When $Y_0=0$, the statistic $\fix^+$ (resp. $\pix^+$) plays no
role. We can then calculate generating functions for signed
(resp. plain) derangements or desarrangements, as shown in the
first two and last rows. The initial polynomial, together with its
two $q$-analogs are reproduced in boldface.

\bigskip
\centerline{\bf 2. Proof of Theorem 1.1}

\medskip
As can be found in ([Co70], p.~9--12)), the generating function for
the derangement numbers~$d_n$ $(n\ge 0)$ is given by
$$\sum_{n\ge 0} d_n{u^n\over n!}=(1-u)^{-1}e^{-u}.
\leqno(2.1)$$ 
An easy
calculation then shows that the polynomials $B_n(Y_0,Y_1,Z)$,
introduced in (1.2), can also be defined by the identity
$$B_n(Y_0,Y_1,Z)=
\sum_{i+j+k+l=n}
{n\choose i,j,k,l}Y_0^i\, Y_1^j\, Z^{j+k}\,d_{k+l}
\quad(n\ge0).\leqno(2.2)
$$

For each signed permutation
$w=x_1x_2\cdots x_n$ let~$A:=\Fix^+w$, $B:=\Fix^-w$,
$C:=\Neg w\setminus \Fix^-w$,
$D:=[n]\setminus(A\cup\overline B\cup \overline C)$. Then
$(A,\overline B,\overline C,D)$ is a sequence of disjoint subsets
of integers, whose union is the interval $[n]:=\{1,2,\ldots,n\}$.
Also the mapping~$\tau$ defined by
$\tau(\overline j)= x_j$ if $\overline j\in C$ and $\tau(j)=x_j$ if
$j\in D$ is a {\it derangement} of the set $C+D$. Hence, $w$ is
completely characterized by the sequence $(A,B,C,D,\tau)$. The
generating polynomial for~$B_n$ by the statistic
$(\fix^+,\fix^-,\negg)$ is then equal to the right-hand side of
(2.2). This proves the first identity of Theorem~1.1.

\medskip
Each signed permutation $w=x_1x_2\cdots x_n$ can be
characterized, either by the four-term sequence
$(\Fix^+w,\Fix^-w,\Neg w,\tau)$, as just described, or
by $(\Pix^+w,\Pix^-w,\Neg w,w^d)$, where
$w^d$ is the desarrangement occurring as the third factor in its
pixed factorization. To construct a bijection~$\phi$ of~$B_n$ onto
$B_n$ such that
$(\fix^-,\fix^+,\negg)\,w=(\pix^-,\pix^+,\negg)\,\phi(w)$
and accordingly prove the second identity of Theorem~1.1,
we only need a bijection $\tau\mapsto f(\tau)$, that maps
each derangement~$\tau$ onto a desarrangement~$f(\tau)$ by
rearranging the letters of~$\tau$. But such a
bijection already exists. It is due to D\'esarm\'enien
({\it op. cit.}). We describe it by means of an example.

\vskip-6pt
Start with a derangement
$\tau=\petitematrice{1&2&3&4&5&6&7&8&9\cr
9&7&4&3&8&2&6&5&1\cr}$ and express it as a product of its
disjoint cycles:
$\tau=(1\,9)(2\,7\,6)(3\,4)(5\,8)$. In each cycle, write the
minimum in {\it second} position:
$\tau=(9\,1)(6\,2\,7)(4\,3)(8\,5)$. Then, reorder the cycles in
such a way that the sequence of those minima, when reading
from left to right, is {\it decreasing}:
$\tau=(8\,5)(4\,3)(6\,2\,7)(9\,1)$. The desarrangement~$f(\tau)$
is derived from the latter expression
by removing the parentheses:
$f(\tau)=8\,5\,4\,3\,6\,2\,7\,9\,1$.

\goodbreak
Let $(\Fix^+w,\Fix^-w,\Neg w,\tau)$ be the sequence
associated with the signed permutation~$w$ and let
$v^-$ (resp.
$v^-$) be the {\it increasing} sequence of the elements
of~$\Fix^-w$ (resp. of~$\Fix^+w$). Then, $v^-\mid
v^+\mid f(\tau)$ is the pixed factorization of
$v^-v^+f(\tau)$ and we may define
$\phi(w)$ by
$$
\phi(w):=v^-v^+f(\tau).\leqno(2.3)
$$
This defines a bijection of $B_n$ onto itself, which has
the further property:
$$
(\Fix^-,\Fix^+,\Neg)\,w=(\Pix^-,\Pix^+,\Neg)\,\phi(w).
\leqno(2.4)
$$ 

For instance, with
\smash{$w=\petitematrice{1&2&3&4&5&6&7&8&9\cr
3&\overline2&8&4&5&\overline1&9&\overline6&7\cr}$} we have
$v^+=4\,5$, $v^-=\overline2$, 
$\tau=\petitematrice{\overline1&3&\overline6&7&8&9\cr
3&8&\overline1&9&\overline6&7\cr}=(9\,7)
(8\,\overline 6\,\overline 1\,3)$ and $f(\tau)=
9\,7\,8\,\overline 6\,\overline 1\,3$. Hence, the pixed
factorization of $\phi(w)$ reads
$\overline 2\mid 4\,5\mid 9\,7\,8\,\overline
6\,\overline 1\,3$ and
$\phi(w)=\overline 2\, 4\,5\, 9\,7\,8\,\overline
6\,\overline 1\,3$.

\bigskip
\centerline{\bf 3. Proof of Theorem 1.2}

\medskip
The length function ``$\ell$" for $B_n$ is expressed in many
ways. We shall use the
following expression derived by Brenti [Br94]. Let
$w=x_1x_2\cdots x_n$ be a signed permutation; its {\it
length}
$\ell(w)$ is defined by
$$\leqalignno{
\ell(w)&:=\inv w+\sum_i|x_i|\,\chi(\, x_i<0),&(3.1)\cr
\noalign{\hbox{where ``inv" designates the usual {\it number of
inversions} for words:}}
\inv w&:=\sum_{1\le i<j\le n}\chi(x_i>x_j).\cr}
$$
The generating polynomial for~$K_n$ (as defined in (1.13))
by ``inv" (resp. for~$D_n$ by ``maj") is denoted by
$K_n(q)$ (resp. $D_n(q)$). 
As was proved in [DeWa93] we 
have:
$$\leqalignno{
K_n(q)&=D_n(q).&(3.2)\cr
\noalign{\hbox{Also}}
\sum_{n\ge 0}{u^n\over (q;q)_n}D_n(q)
&=\Bigl(1-{u\over 1-q}\Bigr)^{-1}\times
(u;q)_\infty,&(3.3)\cr}
$$
as shown by  Wachs [Wa90] in an equivalent form. Another
expression for~$D_n(q)$ will be derived in
section~6, Proposition~6.2.

If $A$ is a finite set of positive integers, let $\tot A$ denote
the sum $\sum a$ $(a\in A)$. For the proof of Theorem~1.2 we
make use of the following classical result, namely that
$q^{N(N+1)/2}{n\brack N}_q$ is equal to the sum
$\sum q^{\tot A}$, where the sum is over all
subsets~$A$ of cardinality~$N$ of the set~$[n]$.

Remember that each signed permutation $w=x_1x_2\ldots x_n$
is characterized by a sequence $(A,B,C,D,\tau)$, where
$A=\Pix^+w$, $B=\Pix^-w$, $C=\Neg w\setminus B$,
$D=[n]\setminus (A\cup\overline B\cup\overline C)$ and
$\tau$ is a {\it desarrangement} of the set $C+D$. Let
$\inv(B,C)$ be the number of pairs of integers
$(i,j)$ such that $i\in B$, $j\in C$ and $i>j$. As
$\inv(B,C)=\inv(\overline C,\overline B)$, we have
$\inv w=\inv(\overline B,\overline C)+\inv(A,D)+\#A\times
\#C+\inv\tau$. From (3.1) it follows that
$$\leqalignno{\noalign{\vskip-6pt}
\ell(w)&=\inv w+\sum_{x_i<0}|x_i|
=\inv w+\tot \overline B+\tot \overline C\cr
\noalign{\vskip-2pt}
&=\tot \overline B+\tot \overline C+\inv(\overline C,\overline
B)+\inv(A,D)+\#A\times
\#C+\inv\tau.\cr}
$$

Denote the right-hand side of (1.10) by $G_n:=G_n(q,Y_0,Y_1,Z)$.
We will calculate 
$G_n(q,Y_0,Y_1,Z)$  by first summing over all sequences
$(A,B,C,D,\tau)$ such that
$\#A=i$, $\#B=j$, $\#C=k$, $\#D=l$. Accordingly, $\tau$ is
a desarrangement of a set of cardinality $k+l$. 
We may write:
$$\leqalignno{\noalign{\vskip-10pt}
G_n&=\kern-8pt\sum_{i+j+k+l=n}
\sum_{(A,B,C,D)}\kern-12pt q^{\tot \overline B+\tot
\overline C+\inv(\overline C,\overline B)
+\inv(A,D)+i\cdot k}\cr
\noalign{\vskip-12pt}
&\kern5.5cm {}\times Y_0^{i}Y_1^{j}Z^{j+k}
\sum_{\tau\in K_{k+l}}q^{\inv\tau}\cr
&=\sum_{m+p=n}\sum_{\scriptstyle
j+k=m,\atop \scriptstyle i+l=p}
\sum_{\scriptstyle \#E=m,\atop
\scriptstyle F=[n]\setminus E}
\sum_{\scriptstyle \overline C+ \overline B=E,\atop
\scriptstyle A+ D=F}\kern-5pt q^{\tot E+\inv
(\overline C,\overline B)+\inv(A,D)+i\cdot k} \cr
\noalign{\vskip-20pt}
&\kern6cm {}\times Y_0^{i}(Y_1Z)^{j}Z^{k}D_{k+l}(q)\cr
\noalign{\smallskip}
&=\sum_{m+p=n}\sum_{\scriptstyle
j+k=m,\atop \scriptstyle i+l=p}
Y_0^i(Y_1Z)^j(Zq^i)^kD_{k+l}(q)\cr
\noalign{\vskip-20pt}
&\kern3.8cm{}\times\kern-8pt\sum_{\scriptstyle \#E=m,\atop
\scriptstyle F=[n]\setminus E}\kern-8ptq^{\tot E}\kern-8pt
\sum_{\scriptstyle \overline C+ \overline B=E,\atop
\scriptstyle A+ D=F}\kern-5pt
q^{\inv(\overline C,\overline B)+\inv(A,D)} \cr
&=\sum_{m+p=n}\sum_{\scriptstyle
j+k=m,\atop \scriptstyle i+l=p}
Y_0^i(Y_1Z)^j(Zq^i)^kD_{k+l}(q)\cr
\noalign{\vskip-20pt}
&\kern3.8cm{}\times q^{m(m+1)/2}{n\brack m}_q\,
{m\brack j,k}_q\,{p\brack i,l}_q.\cr
}
$$
Thus
$$
G_n
= \sum_{i+j+k+l=n}{n\brack i,j,k,l}_qq^{j+k+1\choose 2}
Y_0^i(Y_1Z)^j(Zq^i)^kD_{k+l}(q).\leqno{(3.4)}
$$
Now form the factorial generating function
$$
G(q,Y_0,Y_1,Z; u)
:=\sum_{n\ge 0}{u^n\over (-Zq;q)_n\,(q;q)_n}
G_n(q,Y_0,Y_1,Z).
$$

\goodbreak\noindent
It follows from (3.4) that
$$
\leqalignno{
G(q,Y_0,Y_1,Z;u)&=\sum_{n\ge 0}{1\over (-Zq;q)_n}
\sum_{{i+j+k+l=n}}q^{j+k+1\choose 2}
{(uY_0)^i\over (q;q)_i} {(uY_1Z)^j\over (q;q)_j}\cr
&\kern4.5cm{}\times
u^{n-i-j}{D_{k+l}(q)(Zq^i)^k\over (q;q)_k(q;q)_l}.
\cr}
$$
But ${j+k+1\choose 2}={j+1\choose 2}+(j+1)k+{k\choose 2}$.
Hence
$$
\leqalignno{G(q,Y_0,Y_1,Z;u)
&=\sum_{n\ge 0}{1\over (-Zq;q)_n}
\sum_{m=0}^n\sum_{i+j=m}q^{j+1\choose 2}
{(uY_0)^i\over (q;q)_i}{(uY_1Z)^j\over (q;q)_j}\cr
&\times{u^{n-m}\over
(q;q)_{n-m}}D_{n-m}(q)
 \sum_{k+l=n-m}{n-m\brack k,l}_q
(Zq^{m+1})^kq^{k\choose 2}.\cr
}
$$
Now
$$
\leqalignno{\noalign{\vskip-10pt}
(-Zq^{m+1};q)_{n-m}&=\sum_{k+l=n-m}{n-m\brack k,l}_q
(Zq^{m+1})^kq^{k\choose 2};\cr
\noalign{\vskip-4pt}
\noalign{\hbox{and}}
\noalign{\vskip-4pt}
(-Zq;q)_n&=(-Zq;q)_{m}\, (-Zq^{m+1};q)_{n-m}.\cr
\noalign{\vskip-4pt}}
$$
Hence
$$\leqalignno{\noalign{\vskip-4pt}
G(q,Y_0,Y_1,Z;u)
&=\sum_{n\ge 0}
\sum_{m=0}^n{1\over (-Zq;q)_m}
\sum_{i+j=m} {(uY_0)^i\over (q;q)_i}\,
q^{j+1\choose 2}{(uY_1Z)^j\over (q;q)_j}\cr
\noalign{\vskip-6pt}
&\kern4.5cm{}\times
{u^{n-m}\over (q;q)_{n-m}}D_{n-m}(q)\cr
&=\Bigl(\sum_{n\ge 0}{a_n\,u^n\over (-Zq;q)_n\,(q;q)_n}\Bigr)
\Bigl(\sum_{n\ge 0}{u^n\over (q;q)_n}D_n(q)\Bigr),\cr
\noalign{\hbox{with}}
a_n&=\sum_{i+j=n}{n\brack i,j}_qY_0^iq^{j\choose 2}(qY_1Z)^j
\cr
&=Y_0^n
\sum_{i+j=n}{n\brack i,j}_q(qY_0^{-1}Y_1Z)^j
q^{j\choose 2}\cr
&=Y_0^n(-qY_0^{-1}Y_1Z;q)_n.\cr
}
$$
By taking (3.3) into account this shows that $G(q,Y_0,Y_1,Z;u)$ is equal to
the right-hand side of (1.7) and then 
$G_n(q,Y_0,Y_1,Z)={}^{\ell\kern-1.5pt}B_n(q,Y_0,Y_1,Z)$ holds for
every $n\ge 0$. The proof of Theorem~1.2 is completed. By (3.4)
we also conclude that the identity
$$
(3.5)\ {}^{\ell\kern-1.5pt}B_n(q,Y_0,Y_1,Z)
=\kern-14pt \sum_{i+j+k+l=n}{n\brack i,j,k,l}_qq^{j+k+1\choose 2}
Y_0^i(Y_1Z)^j(Zq^i)^kD_{k+l}(q)
$$
is equivalent to (1.7). As its right-hand side tends to the
right-hand side of (2.2) when $q\rightarrow 1$, we can then
assert that (1.7) specializes into (1.2) for $q=1$.

\bigskip
\centerline{\bf 4. Weighted signed permutations}

\medskip
We use the following notations: if~$c=c_1c_2\cdots c_n$
is a word, whose letters are nonnegative integers, let 
$\lambda(c):=n$ be the {\it length} of~$c$,
$\tot c:=c_1+c_2+\cdots +c_n$ the {\it sum} of its letters and
$\odd c$ the number of its {\it odd} letters. 
Furthermore, $\NIW_n$ (resp. $\NIW_n(s)$) designates the set
of all {\it nonincreasing} words of length~$n$, whose letters
are nonnegative integers (resp. nonnegative integers at most
equal to~$s$). Also let $\NIW^e_n(s)$ (resp. $\DW_n^o(s)$) be
the subset of $\NIW_n(s)$  of the nonincreasing (resp. strictly
decreasing) words all letters of which are {\it even} (resp. {\it
odd\/}). 

\medskip
Next, each pair
${c\choose w}$ is called a {\it weighted signed permutation}
of order~$n$ if  the four properties 
({\sl wsp1})--({\sl wsp4}) hold:

({\sl wsp1}) $c$ is a word $c_1c_2\cdots c_n$ from $\NIW_n$; 

({\sl wsp2}) $w$ is a signed permutation $x_1x_2\cdots x_n$ from
$B_n$;

({\sl wsp3}) $c_k=c_{k+1}\Rightarrow x_k<x_{k+1}$ for all
$k=1,2,\ldots, n-1$;

({\sl wsp4}) $x_k$ is positive (resp. negative) whenever $c_k$ is
even (resp. odd).

\medskip
When $w$ has no fixed points, either negative or positive, we
say that $c\choose w$ is a {\it weighted signed derangement}.
The set of weighted signed permutations (resp. derangements)
${c\choose w}={c_1c_2\cdots c_n\choose x_1x_2\cdots x_n}$
of order~$n$ is denoted by $\WSP_n$ (resp. by~$\WSD_n$).
The subset of all those weighted signed permutations (resp.
derangements) such that
$c_1\le s$ is denoted by $\WSP_n(s)$ (resp. by~$\WSD_n(s)$).

\medskip
For example, the following pair 
$$
\matrice{\cr
\displaystyle{c\choose w}=\cr}
\left(\matrice{\strut 1&2&\vrule\ 3&\vrule\ 4&5&6&\vrule\ 7
&8&9&\vrule\ 10&\vrule\ 11&12&\vrule\ 13\cr
\strut10&10&\vrule\ 9&\vrule\ 7&7&7&\vrule\ 4&4&4
&\vrule\ \hphantom{0} 3&
\vrule\ \hphantom{0} 2&\hphantom{0} 2&\vrule\ \hphantom{0} 1\cr
\strut 1&2&\vrule\ \overline 7&\vrule\ 
\overline6&\overline 5&\overline 4&\vrule\ 3&8&9&\vrule\
\overline{10}&\vrule\ 12&13&\vrule\ 
\overline {11}\cr}\right)
$$
is a weighted signed permutation of order~13. It has four
positive fixed points (1, 2, 8, 9) and two negative fixed points
($\overline 5$, $\overline{10}$).

\proclaim Proposition 4.1. With each weighted signed
permutation ${c\choose w}$ from the set $\WSP_n(s)$
can be associated a unique sequence \smash{$(i,j,k,{c'\choose
w'},v^e,v^o)$} such that\hfil\break
\indent{\rm (1)} $i$, $j$, $k$ are nonnegative integers of
sum~$n$;\hfil\break
\indent{\rm (2)} ${c'\choose w'}$ is a weighted signed
derangement from the set $\WSD_i(s)$;\hfil\break
\indent{\rm (3)} $v^e$ is a nonincreasing word with even
letters from the set $\NIW_j^e(s)$;\hfil\break
\indent{\rm (4)} $v^o$ is a decreasing word with odd
letters from the set $\DW_k^o(s)$;\hfil\break
having the following properties:
$$
\eqalign{
\tot c=\tot c'+\tot v^e&+\tot v^o;\quad
\negg w=\negg w'+\lambda (v^o);\cr
\fix w^+=\lambda (v^e);&\quad
\fix^-w=\lambda (v^o).\cr}\leqno(4.1)
$$

\goodbreak
The bijection ${c\choose w}\mapsto ({c'\choose w'},v^e,v^o)$
is quite natural to define. Only its reverse requires some
attention. To get the latter three-term sequence from
$c\choose w$ proceed as follows: 

(a) let ${l_1}$, \dots~, ${l_{\alpha}}$ (resp. ${m_1}$, \dots~,
${m_{\beta}}$) be the increasing sequence of the integers~$l_i$
(resp. $m_i$) such that $x_{l_i}$ (resp. $x_{m_i}$) is a positive
(resp. negative) fixed point of~$w$;

(b) define: $v^e:=c_{l_1}\cdots c_{l_{\alpha}}$ and
$v^o:=c_{m_1}\cdots c_{m_{\beta}}$;

(c) remove all the columns \smash{${c_{l_1}\choose
x_{l_1}},\ldots,{c_{l_{\alpha}}\choose
x_{l_{\alpha}}}$, ${c_{m_1}\choose
x_{m_1}},\ldots,{c_{m_{\beta}}\choose x_{m_{\beta}}}$} from
$c\choose w$ and let $c'$ be the nonincreasing word derived
from~$c$ after the removal;

(d) once the letters $x_{l_1},\ldots, x_{l_{\alpha}}, 
x_{m_1},\ldots, x_{m_{\beta}}$ have been removed from the
signed permutation~$w$ the remaining ones form a signed
permutation of a subset~$A$ of~$[n]$, of cardinality
$n-\alpha-\beta$. Using the unique increasing bijection~$\phi$
of~$A$ onto the interval $[n-\alpha-\beta]$ replace each
remaining letter~$x_i$ by $\phi(x_i)$ if $x_i>0$ or by
$-\phi(-x_i)$ if $x_i<0$. Let~$w'$ be the signed derangement of
order $n-\alpha-\beta$ thereby obtained.

\medskip
For instance, with the above weighted signed permutation we have:
$v^e=10, 10, 4, 4$ and $v^o=7, 3$. After removing the fixed
point columns we obtain:
$$
\left(\matrice{\strut 3&\vrule\ 4&6&\vrule\ 7
&\vrule\ 11&12&\vrule\ 13\cr
\strut9&\vrule\ 7&7&\vrule\ 4
&\vrule\ \hphantom0 2&\hphantom0 2&\vrule\ \hphantom0 1\cr
\strut \overline 7&\vrule\ \overline 6&\overline 4&\vrule\
3&\vrule\ 12&13&\vrule\ 
\overline {11}\cr}\right)
\hbox{ and then }
\matrice{\cr
\displaystyle{c'\choose w'}=\cr}
\left(\matrice{\strut 1&\vrule\ 2&3&\vrule\ 4
&\vrule\ 5&6&\vrule\ 7\cr
\strut9&\vrule\ 7&7&\vrule\ 4
&\vrule\ 2&2&\vrule\ 1\cr
\strut \overline 4&\vrule\ \overline 3&\overline 2&\vrule\
1&\vrule\ 6&7&\vrule\ 
\overline {5}\cr}\right).
$$

There is no difficulty verifying that the properties listed in
(4.1) hold.  For reconstructing $c\choose w$ from the sequence
$({c'\choose w'},v^e,v^o)$ consider the nonincreasing
rearrangement of the juxtaposition product $v^ev^o$ in the
form $b_1^{h_1}\cdots b_m^{h_m}$, where $b_1>\cdots >b_m$ and
$h_i\ge 1$ (resp. $h_i=1$) if $b_i$ is even (resp. odd). The
pair $c'\choose w'$ being decomposed into matrix blocks, as
shown in the example, each letter~$b_i$ indicates where
the~$h_i$ fixed point columns are to be inserted. We do not
give more details and simply illustrate the construction with
the running example.

\medskip
With the previous example $b_1^{h_1}\cdots b_m^{h_m}
=10^2\,7\,4^2\,3$. First, implement $10^2$:
$$\displaylines{\noalign{\vskip-10pt}
\left(\matrice{\strut 1&2&\vrule\ 3&\vrule\ 4&5&\vrule\ 6
&\vrule\ 7&8&\vrule\ 9\cr
\strut\bf10&\bf10&\vrule\ 9&\vrule\ 7&7&\vrule\ 4
&\vrule\ 2&2&\vrule\ 1\cr
\strut 1&2&\vrule\ \overline 6&\vrule\ \overline 5&\overline
4&\vrule\ 3&\vrule\ 8&9&\vrule\ 
\overline {7}\cr}\right);\cr
\noalign{\hbox{then 7:}}
\left(\matrice{\strut 1&2&\vrule\ 3&\vrule\ 4&5&6&\vrule\ 7
&\vrule\ 8&9&\vrule\ 10\cr
\strut10&10&\vrule\ 9&\vrule\ 7&\bf7&7&\vrule\ 4
&\vrule\ 2&2&\vrule\ \hphantom0 1\cr
\strut 1&2&\vrule\ \overline 7&\vrule\ \overline 6&\overline
5&\overline 4&\vrule\ 3&\vrule\ 9&10&\vrule\ \hphantom0 
\overline {8}\cr}\right);\cr}
$$
notice that because of condition ({\sl wsp3}) the letter {\bf 7} is
to be inserted in {\it second} position in the third block;

\goodbreak
\noindent then insert $4^2$:
$$
\displaylines{\noalign{\vskip-4pt}
\left(\matrice{\strut 1&2&\vrule\ 3&\vrule\ 4&5&6&\vrule\ 7
&8&9&\vrule\ 10&11&\vrule\ 12\cr
\strut10&10&\vrule\ 9&\vrule\ 7&7&7&\vrule\ 4&\bf4&\bf4
&\vrule\ \hphantom0 2&2&\vrule\ \hphantom0 1\cr
\strut 1&2&\vrule\ \overline 7&\vrule\ \overline 6&\overline
5&\overline 4&\vrule\ 3&8&9&\vrule\ 11&12&\vrule\ 
\overline {10}\cr}\right).\cr
}
$$
The implementation of 3 gives back the original weighted
signed permutation~$c\choose w$.

\goodbreak
\bigskip
\centerline{\bf 5. Proof of Theorem 1.3}

\medskip
It is $q$-routine (see, e.g., [An76, chap.~3]) to prove the
following identities, where~$v_1$ is the first letter of~$v$:
$$
\displaylines{
{1\over (u;q)_N}=\sum_{n\ge 0} {N+n-1\brack
n}_q\,u^n;\qquad
{N+n\brack n}_q=\sum_{v\in \hbox{\sevenrm NIW}_n(N)}q^{\tot
v};\cr 
{1\over (u;q)_{N+1}}=\sum_{n\ge 0}u^n
\sum_{v\in \hbox{\sevenrm NIW}_n(N)}q^{\tot v}
={1\over 1-u}\sum_{v\in \hbox{\sevenrm NIW}_n}q^{\tot
v}u^{v_1};\cr 
\rlap{(5.1)}\hfill {1\over
(u;q^2)_{\lfloor s/2\rfloor+1}} =\sum_{n\ge 0}u^n
\sum_{v^e\in \hbox{\sevenrm NIW}_n^e(s)}
q^{\tot v^e};\hfill\cr
\rlap{(5.2)}\hfill
(-uq;q^2)_{\lfloor (s+1)/2\rfloor}
=\sum_{n\ge 0}u^n
\sum_{v^o\in \hbox{\sevenrm DW}_n^o(s)}
q^{\tot v^o}.\hfill\cr
}
$$

The last two formulas and Proposition 4.1 are now used to
calculate the generating function for the weighted signed
permutations. The symbols $\NIW^e(s)$, $\DW^o(s)$, $\WSP(s)$,
$\WSD(s)$ designate the unions for $n\ge 0$ of the
corresponding symbols with an $n$-subscript. 

\proclaim Proposition 5.2. The following identity holds:
$$
\displaylines{
(5.3)\quad 
\sum_{n\ge 0}u^n
\kern-15pt\sum_{\quad{c\choose w}\in\hbox{\sevenrm
WSP}_n(s)}
\kern-18pt
q^{\tot c}Y_0^{\fix^+ w}Y_1^{\fix^-w}Z^{\negg
w}\hfill\cr
\hfill{}
=
{(u;q^2)_{\lfloor s/2\rfloor+1}\over (uY_0;q^2)_{\lfloor
s/2\rfloor+1}}{(-uqY_1Z;q^2)_{\lfloor (s+1)/2\rfloor}\over
(-uqZ;q^2)_{\lfloor (s+1)/2\rfloor}} \times
\sum_{n\ge 0}u^n
\kern-18pt\sum_{\quad{c\choose w}\in\hbox{\sevenrm
WSP}_n(s)}
\kern-24pt
q^{\tot c}Z^{\negg w}\
.\cr }
$$

{\it Proof}.\quad
First, summing over $(w^e,w^o,{c\choose w})\in
\NIW(s)\times\DW(s)\times \WSP(s)$, we have
$$\displaylines{
\sum_{w^e,w^o,{c\choose w}}
u^{\lambda (w^e)}
q^{\tot w^e}\times (uZ)^{\lambda (w^o)}q^{\tot w^o}
\times
u^{\lambda (c)}
q^{\tot c}Y_0^{\fix^+w}
Y_1^{\fix^-w}Z^{\negg w}\hfill\cr
\noalign{\vskip-6pt}
(5.4)\hfill{}={(-uqZ;q^2)_{\lfloor (s+1)/2\rfloor}
\over (u;q^2)_{\lfloor s/2\rfloor+1}}\times
\sum_{{c\choose w}}
u^{\lambda (c)}
q^{\tot c}Y_0^{\fix^+w}
Y_1^{\fix^-w}Z^{\negg w}\cr
}
$$
by (5.1) and (5.2). Now,
Proposition 4.1 implies that the initial expression can also be
summed over five-term sequences $({c'\choose w'},v^e,v^o, w^e,
w^o)$ from $\WSD(s)\times \NIW^e(s)\times \DW^o(s)\times
\NIW^e(s)
\times \DW^o(s)$ in the form
$$
\displaylines{\kern-30pt
\sum_{\qquad\quad{c'\choose w'},v^e,v^o, w^e,
w^o}\kern-27pt u^{\lambda( c')}
q^{\tot c'}Z^{\negg w'}\times
(uY_0)^{\lambda (v^e)}q^{\tot v^e}\times
(uY_1Z)^{\lambda (v^o)}q^{\tot v^o}\hfill\cr
\noalign{\vskip-12pt}
\hfill{}\times u^{\lambda (w^e)}q^{\tot w^e} 
\times (uZ)^{\lambda
(w^o)}q^{\tot w^o}\cr
\kern.5cm{}=\sum_{v^e,v^o}
(uY_0)^{\lambda (v^e)}q^{\tot v^e}\times
(uY_1Z)^{\lambda (v^o)}q^{\tot v^o}
\hfill\cr
\noalign{\vskip-6pt}
\hfill{}\times
\sum_{{c'\choose w'}, w^e,
w^o} u^{\lambda (c')}
q^{\tot c'}Z^{\negg w'}
\times u^{\lambda (w^e)}q^{\tot w^e} 
\times (uZ)^{\lambda
(w^o)}q^{\tot w^o}.\cr
}
$$
The first summation can be evaluated by (5.1) and (5.2), while
by Proposition 4.1 again the second sum can be expressed
as a sum over weighted signed permutations ${c\choose w}\in
\WSP(s)$. Therefore, the initial sum is also equal to
$$
{(-uqY_1Z;q^2)_{\lfloor (s+1)/2\rfloor}\over
(uY_0;q^2)_{\lfloor
s/2\rfloor+1}}\times\kern-16pt
\sum_{\quad{c\choose w}\in\hbox{\sevenrm
WSP}(s)}
u^{\lambda (c)}q^{\tot c}Z^{\negg w}.
\leqno(5.5)
$$
Identity (5.3) follows by equating (5.4) with (5.5).\cqfd

\proclaim Proposition 5.3. The following identity holds:
$$
\sum_{n\ge 0}u^n
\kern-15pt\sum_{\quad{c\choose w}\in\hbox{\sevenrm
WSP}_n(s)}
\kern-18pt
q^{\tot c}Z^{\negg w}
=\Bigl( 1-u\sum_{i=0}^s q^iZ\strut^{\chi (i\ {\rm
odd})}\Bigr)^{-1}.\leqno(5.6)
$$

\proof
For proving the equivalent identity
$$
\sum_{\quad{c\choose w}\in\hbox{\sevenrm
WSP}_n(s)}\kern-18pt q^{\tot c}Z^{\negg w}
=\Bigl(\sum_{i=0}^sq^iZ\strut^{\chi (i\ {\rm odd})}
\Bigr)^n\quad (n\ge 0)\leqno(5.7)
$$
it suffices to  construct a bijection ${c\choose w}\mapsto
d$ of
$\WSP_n(s)$ onto $\{0,1,\ldots,s\}^n$ such that
$\tot c=\tot d$~and $\negg w=\odd d$. This bijection is
one of the main ingredients of the {\it MacMahon Verfahren} for
signed permutations that has been  fully described in [FoHa05,
\S\kern2pt 4]. We simply recall the construction of the
bijection by means of an example.

Start with
$\petitematrice{c\cr w\cr}=
\petitematrice{
10&9&7&4&4&2&2&1&1\cr
1&\overline4&\overline{3}&2&5&6&8&
\overline{9}&\overline7\cr}$.
Then, form the two-matrix
$\petitematrice{
10&9&7&4&4&2&2&1&1\cr
1&4&{3}&2&5&6&8&
{9}&7\cr}$, where the {\it negative} integers on the bottom row
have been replaced by their opposite values. Next, rearrange
its columns in such a way that the bottom row is precisely
$1\,2\,\ldots\,n$. The word~$d$ is defined to be the top row
in the resulting matrix. Here
${\displaystyle{d\choose {\rm Id}}}\!\!
=\!\!\petitematrice{\!1\!0&4&7&9&4&2&1&2&1\cr
1&2&3&4&5&6&7&8&9\cr}$. 
As~$d$ is a rearrangement of~$c$,
we have $\tot c=\tot d$ and $\negg w=\odd d$. 
For reconstructing the pair $c\choose w$ from
$d=d_1d_2\cdots d_n$ simply make a full use of condition
({\sl wsp3}).

Using the properties of this bijection we have:
$$\displaylines{\noalign{\vskip-7pt}
\sum_{\kern-5pt \quad{c\choose w}\in\hbox{\sevenrm
WSP}_n(s)}\kern-18pt q^{\tot c}Z^{\negg w}
=\kern-15pt \sum_{d\in \{0,1,\ldots,s\}^n}\kern-12pt
q^{\tot d}Z^{\odd d}\
=\kern-12pt\sum_{d\in  \{0,1,\ldots,s\}^n}
\prod_{i=1}^n q^{d_i}
Z\strut^{\chi (d_i\ {\rm odd})}\hfill\cr
\noalign{\vskip-3pt}
\hfill{}=\prod_{i=1}^n\sum_{d_i\in \{0,1,\ldots,s\}}
 \kern-15pt q^{d_i}
Z^{\chi (d_i\ {\rm odd})}
=\Bigl(\sum_{i=0}^sq^iZ\strut^{\chi (i\ {\rm odd})}
\Bigr)^n\kern-5pt .\qed\quad\cr}
$$

Let $G_n:=G_n(t,q,Y_0,Y_1,Z)$ denote the right-hand side of (1.12)
in the statement of Theorem~1.3.

\proclaim Proposition 5.4. 
Let $G_n\!:=G_n(t,q,Y_0,Y_1,Z)$ denote the right-hand~side of
$(1.12)$ in the statement of Theorem~$1.3$. Then
$$
{1+t\over
(t^2;q^2)_{n+1}}\,G_n
=\sum_{s\ge 0} t^s\sum_{{c\choose w}\in
\hbox{\sevenrm WSP}_n(s)}
\kern-10pt q^{\tot c}Y_0^{\fix^+ w}
Y_1^{\fix^- w}Z^{\negg
w}.\leqno(5.8)
$$

\proof
A very similar calculation has been made in the proof of
Theorem~4.1 in [FoHa05]. We also make use of the identities
on the $q$-ascending factorials that were recalled in the
beginning of this section. First,
$$\eqalignno{\noalign{\vskip-6pt}
{1+t\over(t^2;q^2)_{n+1}}
&=\sum_{r'\ge 0}(t^{2r'}+t^{2r'+1}){n+r'\brack r'}_{q^2}\cr
&=\sum_{r\ge 0}t^r{n+\lfloor r/2\rfloor \brack
\lfloor r/2\rfloor}_{q^2}
=\sum_{r\ge 0}t^r
\sum_{b\in \hbox{\sixrm NIW}_n(\lfloor
r/2\rfloor)}\kern-12pt q^{2\tot b}.\cr
\noalign{\vskip-6pt}}
$$
Then,
$$
\leqalignno{\noalign{\vskip-6pt}
{1+t\over
(t^2;q^2)_{n+1}}\,G_n
&=\sum_{r\ge 0}t^r
\sum_{\scriptstyle b\in \hbox{\sixrm NIW}_n,\atop
\scriptstyle 2b_1\le
r}\kern-12pt q^{2\tot b}
\sum_{w\in B_n} t^{\fdes w}q^{\fmaj w}
Y_0^{\fix^+ w}Y_1^{\fix^-w}Z^{\negg w}\hfill\cr
&=
\sum_{s\ge 0}t^s
\sum_{\scriptstyle b\in \hbox{\sixrm NIW}_n,\,w\in B_n
\atop \scriptstyle
2b_1+\fdes w\le s}
\kern-12pt q^{2\tot b+\fmaj w}
Y_0^{\fix^+ w}Y_1^{\fix^-w}Z^{\negg w}.\hfill
\cr
}
$$
As proved in [FoHa05, \S\kern2pt 4] to each ${c\choose w}
={c_1\cdots c_n\choose x_1\cdots x_n}\in
\WSP_n(s)$ there corresponds a unique $b=b_1\cdots b_n\in
\NIW_n$ such that $2b_1+\fdes w=c_1$ and $2\tot b+\fmaj
w=\tot c$. Moreover, the mapping ${c\choose w}\mapsto
(b,w)$ is a bijection of $\WSP_n(s)$ onto the set of all pairs 
$(b,w)$ such that $b=b_1\cdots b_n\in
\NIW_n$ and $w\in B_n$ with the property that
$2b_1+\fdes w\le s$.

The word~$b$ is determined as follows: write the signed
permutation~$w$ as a linear word
$w=x_1x_2\ldots x_n$ and for
each
$k=1,2,\ldots,n$ let $z_k$ be the number of descents 
$(x_i>x_{i+1}$) in the right
factor $x_kx_{k+1}\cdots x_n$ and let $\epsilon_k$
be equal to~0 or~1 depending on whether $x_k$ is positive
or negative. Also for each $k=1,2,\ldots ,n$ define
$a_k:=(c_k-\epsilon_k)/2$, $b_k:=(a_k-z_k)$ and form the
word $b=b_1\cdots b_n$.

For example, 
$$
\matrice{{\rm Id}&=&1&2&3&4&5&6&7&8&9&10\cr
c&=&9&7&7&4&4&4&2&2&1&1\cr
w&=&
\overline4&\overline3&\overline2&1&5&6&8&9
&\overline{10}&\overline7\cr
z&=&1&1&1&1&1&1&1&1&0&0\cr
\epsilon&=&1&1&1&0&0&0&0&0&1&1\cr
a&=&4&3&3&2&2&2&1&1&0&0\cr
b&=&3&2&2&1&1&1&0&0&0&0\cr}
$$
Pursuing the above calculation we get (5.8).\qed

\bigskip
We can complete the proof of Theorem 1.3:
$$\displaylines{\noalign{\vskip-2pt}
\quad
\sum_{n\ge 0}(1+t)G_n(t,q,Y_0,Y_1,Z){u^n\over
(t^2;q^2)_{n+1}}
\hfill\cr
\noalign{\vskip-4pt}
\kern 1cm{}=\sum_{s\ge 0}t^s\sum_{n\ge 0}u^n
\sum_{{c\choose w}\in
\hbox{\sevenrm WSP}_n(s)}
\kern-10pt q^{\tot c}Y_0^{\fix^+ w}
Y_1^{\fix^- w}Z^{\negg
w}\hfill\hbox{\rm [by (5.8)]}\cr
\noalign{\vskip-4pt}
\kern 1cm{}
=\sum_{s\ge 0}t^s
{(u;q^2)_{\lfloor s/2\rfloor+1}\over (uY_0;q^2)_{\lfloor
s/2\rfloor+1}}{(-uqY_1Z;q^2)_{\lfloor (s+1)/2\rfloor}\over
(-uqZ;q^2)_{\lfloor (s+1)/2\rfloor}} \hfill\cr
\hfill{}\times \sum_{n\ge 0}u^n \kern-10pt
\sum_{\quad{c\choose w}\in\hbox{\sevenrm
WSP}_n(s)}\kern-10pt
q^{\tot c}Z^{\negg w}\quad\hbox{\rm
[by (5.3)]}\cr
\noalign{\vskip-8pt}
\kern 1cm{}
=\sum_{s\ge 0} t^s
\Bigl(1-u\sum_{i=0}^s q^iZ\strut^{\chi (i\ {\rm
odd})}\Bigr)^{-1}\hfill\cr
\noalign{\vskip-2pt}
\hfill{}\times{(u;q^2)_{\lfloor s/2\rfloor+1}\over
(uY_0;q^2)_{\lfloor s/2\rfloor+1}}{(-uqY_1Z;q^2)_{\lfloor
(s+1)/2\rfloor}\over (-uqZ;q^2)_{\lfloor (s+1)/2\rfloor}}
\quad\hbox{\rm [by (5.6)].}\cr }$$
Hence, $G_n(t,q,Y_0,Y_1,Z)=B_n(t,q,Y_0,Y_1,Z)$ for all $n\ge 0$.\qed

\bigskip
\centerline{\bf 6. Specializations}

\medskip
For deriving the specializations of the polynomials 
${}^{\ell\kern-1.5pt}B_n(q,Y_0,Y_1,Z)$ and $B_n(t,q,Y_0,Y_1,Z)$
with their combinatorial interpretations we refer to the diagram
displayed in Fig.~1. Those two polynomials are now regarded as
generating polynomials for~$B_n$ by the multivariable statistics
$(\ell,\pix^+,\pix^-,\negg)$ and 
$(\fdes,\fmaj,\fix^+,\fix^-,\negg)$,  their factorial
generating functions being given by (1.7) and (1.9), respectively. 

First, identity (1.8) is deduced from (1.9) by the traditional token
that consists of multiplying (1.9) by $(1-t)$ and making $t=1$. 
Accordingly, $B_n(q,Y_0,Y_1,Z)$ occurring in (1.8) is the
generating polynomial for the group~$B_n$ by the statistic
$(\fmaj,\fix^+,\fix^-,\negg)$. 

Now, let
$$
B(q,Y_0,Y_1,Z;u):=
\sum_{n\ge 0}{u^n\over
(q^2;q^2)_n}B_n(q,Y_0,Y_1,Z).\leqno(6.1)
$$
The {\it involution} of~$B_n$ defined by
$w=x_1x_2\cdots x_n\mapsto \overline w:=\overline x_1
\overline x_2\cdots \overline x_n$ has the following properties:
$$\leqalignno{
\fmaj w+\fmaj \overline w=n^2;&\quad 
\negg w+\negg\overline w=n;&(6.2)\cr
\fix^+w=\fix^- \overline w;&\quad\fix^- w=\fix^+\overline
w.&(6.3)\cr}$$
Consequently, the duality between positive and negative fixed
points must be reflected in the expression of
$B(q,Y_0,Y_1,Z;u)$ itself, as shown next.

\proclaim Proposition 6.1. We have:
$$B(q,Y_0,Y_1,Z;u)=B(q^{-1},Y_1,Y_0,Z^{-1};-uq^{-1}Z).
\leqno(6.4)
$$

\proof
The combinatorial proof consists of using the relations written
in (6.2), (6.3) and easily derive the identity
$$
B_n(q,Y_0,Y_1,Z)=q^{n^2}Z^n \,B_n(q^{-1},Y_1,Y_0,Z^{-1}).
\leqno(6.5)
$$
With this new expression for the generating polynomial
identity~(6.1) becomes
$$
B(q,Y_0,Y_1,Z;u)=
\sum_{n\ge 0}{(-uq^{-1}Z)^n\over
(q^{-2};q^{-2})_n}B_n(q^{-1},Y_1,Y_0,Z^{-1}),
$$
which implies (6.4).

The analytical proof consists of showing that the right-hand
side of identity (1.8) is invariant under the transformation
$$
(q,Y_0,Y_1,Z,u)\mapsto (q^{-1},Y_1,Y_0,Z^{-1},-uq^{-1}Z).
$$
The factor $1-u(1+qZ)/(1-q^2)$ is clearly invariant.
As for the other two factors it suffices to expand them by
means of the $q$-binomial theorem ([GaRa90], p.~7) and
observe that they are simply permuted when the
transformation is applied.\cqfd

\goodbreak
\medskip
The polynomial $D_n^B(t,q,Y_1,Z):=B_n(t,q,0,Y_1,Z)$ 
(resp. $D_n^B(q,Y_1,Z)\hfil\break:=B_n(q,0,Y_1,Z)$) is the
generating polynomial for the set $D_n^B$ of the {\it signed
derangements} by the statistic $(\fdes,\fmaj,\fix^-,\negg)$
(resp. $(\fmaj,\fix^-,\negg)$).
Their factorial generating functions are obtained by letting
$Y_0=0$ in (1.9) and (1.8), respectively.

Let $Y_0=0$, $Y_1=1$ in (1.8). We then obtain
the factorial generating function for the polynomials
$D_n^B(q,Z):=\sum q^{\fmaj w}Z^{\negg w}$
$(w\in D_n^B)$ in the form
$$
\sum_{n\ge 0}{u^n\over (q^2;q^2)_n}D_n^B(q,Z)=
\Bigl(1-u{1+qZ\over 1-q^2}\Bigr)^{-1}
\times (u;q^2)_\infty.
\leqno(6.6)
$$
It is worth writing the equivalent forms of that identity:
$$
\displaylines{
\rlap{(6.7)}
\kern1cm {(q^2;q^2)_n\over (1-q^2)^n}
\,(1+qZ)^n=\sum_{k=0}^n{\,n\,\brack
k}_{q^2}\kern-5pt D_k^B(q,Z)\quad(n\ge 0);\hfill\cr
\rlap{(6.8)}\kern1cm
D_n^B(q,Z)\!=\!\sum_{k=0}^n {n\brack
k}_{q^2}\kern-9pt (-1)^kq^{k(k-1)}{(q^2;q^2)_{n-k}\over
(1-q^2)^{n-k}} (1+qZ)^{n-k}\quad(n\ge 0);\hfill\cr
\rlap{(6.9)}\kern1cm
D_0^B(q,Z)=1,\quad{\rm and\ for\ }n\ge 0\hfill\cr
\kern1cm D_{n+1}^B(q,Z)=(1+qZ){1-q^{2n+2}\over
1-q^2}D_n^B(q,Z)+(-1)^{n+1}q^{n(n+1)}.\hfill\cr
\rlap{(6.10)}\kern1cm
D_0^B(q,Z)=1,\quad D_1^B(q,Z)=Zq,
\quad{\rm and\ for\ }n\ge 1\hfill\cr
\kern1cm
D_{n+1}^B(q,Z)=
\Bigl({1-q^{2n}\over 1-q^2}+qZ{1-q^{2n+2}\over
1-q^2}\Bigr)D_{n}^B(q,Z)\hfill\cr
\kern5cm{}
+(1+qZ)q^{2n}{1-q^{2n}\over 1-q^2}
D_{n-1}^B(q,Z).
\hfill\cr}
$$
Note that (6.8) is derived from (6.6) by taking the
coefficients of~$u^n$ on both sides. Next, multiply both sides
of (6.6) by the second $q^2$-exponential
$E_{q^2}(-u)$ and look for the coefficients of~$u^n$ on both
sides. This yields (6.7).
Now, write (6.6) in the form
$$E_{q^2}(-u)=
\Bigl(1-u{1+qZ\over 1-q^2}\Bigr)
\sum_{n\ge 0}{u^n\over (q^2;q^2)_n}D_n^B(q,Z)\leqno(6.11)
$$
and take the coefficients of~$u^n$ on both sides. This yields (6.9).
Finally, (6.10) is a simple consequence of (6.9). 

When $Z=1$, formulas (6.6),
(6.8), (6.9) have been proved by Chow [Ch06] with
$D_n^B(q)=\sum_w q^{\fmaj w}$ $(w\in D_n^B$).

Now the polynomial $K_n^B(q,Y_1,Z):={}^{\ell\kern-1.5pt}B_n(q,0,Y_1,Z)$
is the generating polynomial for the set $K_n^B$ of the {\it
signed desarrangements} by the statistic $(\ell,\pix^-,\negg)$.
From (1.7) we get
$$
\displaylines{(6.12)\quad
\sum_{n\ge 0}{u^n\over (-Zq;q)_n\,(q;q)_n}
K_n^B(q,Y_1,Z)\hfill\cr
\noalign{\vskip-6pt}
\hfill{}
=\Bigl(1-{u\over 1-q}\Bigr)^{-1}\times
(u;q)_\infty\,\Bigl(\sum_{n\ge 0}
{q^{n+1\choose 2}(Y_1Zu)^n\over (-Zq;q)_n\,(q;q)_n}\Bigr).
\quad\cr}
$$

When the variable $Z$ is given the zero value, the polynomials
in the second column of Fig.~1 are mapped on generating
polynomials for the {\it symmetric group}, listed in the third
column. Also the variable~$Y_1$ vanishes. Let
$A_n(t,q,Y_0):=B_n(t^{1/2},q^{1/2},Y_0,0,0)$. Then
$$\displaylines{\rlap{(6.13)}\hfill
A_n(t,q,Y_0)=\sum_{\sigma\in
{\goth S}_n} t^{\des \sigma}q^{\maj\sigma}Y_0^{\fix \sigma}
\quad (\fix:=\fix^+).
\hfill\cr
\noalign{\hbox{Identity (1.9)
specializes into}}
\rlap{(6.14)}\hfill\quad
\sum_{n\ge 0}A_n(t,q,Y_0){u^n\over (t;q)_{n+1}}
=\sum_{s\ge 0}t^s\Bigl(
1-u\sum_{i=0}^sq^i\Bigr)^{-1} {(u;q)_{s+1}\over
(uY_0;q)_{s+1}},
\hfill\cr}
$$
an identity derived by Gessel and Reutenauer ([GeRe93],
Theorem~8.4) by means of a quasi-symmetric function
technique. Note that they wrote their formula for
``$1+\des$" and not for~``$\des$."

Multiply (6.14) by $(1-t)$ and let $t:=1$, or let $Z:=0$ and
$q^2$ be replaced by~$q$ in (1.8). Also, let
\smash{$A_n(q,Y_0):=\sum_{\sigma}q^{\maj\sigma} Y^{\fix \sigma}$}
$(\sigma\in {\goth S}_n)$; we get
$$
\sum_{n\ge 0}{u^n\over (q;q)_n}A_n(q,Y_0)
=\Bigl(1-{u\over 1-q}\Bigr)^{-1}{(u;q)_\infty\over
(uY_0;q)_\infty},\leqno(6.15)
$$
an identity derived by Gessel and Reutenauer [GeRe93] and 
also by Clarke {\it et al.} [ClHaZe97] by means of a $q$-Seidel
matrix approach.

We do not write the specialization of (6.14) when $Y_0:=0$ to
obtain the generating function for the polynomials
$D_n(t,q):=\sum\limits_{\sigma\in
D_n}t^{\des\sigma}q^{\maj\sigma}$. As for the polynomial
$D_n(q):=\sum\limits_{\sigma\in D_n}q^{\maj\sigma}$, it has
several analytical expressions, which can all be derived from
(6.7)--(6.10) by letting $Z:=0$ and $q^2$ being replaced by~$q$.
We only write the identity which corresponds to (6.7)
$$\displaylines{\noalign{\vskip-6pt}
(6.16)\qquad
D_0(q)=1\quad
\hbox{and}\quad
{(q;q)_n\over (1-q)^n}=\sum_{k=0}^n{n\brack k}_qD_k(q)
\hbox{\quad for $n\ge 1$,}\quad
\hfill
\cr
\noalign{\hbox{which is then equivalent to the identity}}
(6.17)\kern1.5cm
e_q(u)\sum_{n\ge 0}{u^n\over (q;q)_n}D_n(q)
=\Bigl(1-{u\over 1-q}\Bigr)^{-1}.\hfill\cr}
$$

\goodbreak
The specialization of (6.8) for $Z:=0$ and $q^2$ replaced by~$q$ was
originally proved by Wachs [Wa98] and again recently by Chen and Xu
[ChXu06]. Those two authors make use of the now classical {\it MacMahon
Verfahren}, that has been exploited in several papers and further extended
to the case of signed permutation, as described in our previous paper
[FoHa05].

In the next proposition we show that $D_n(q)$ can be expressed as
a polynomial in~$q$ with {\it positive integral} coefficients. In
the same manner, the usual derangement number~$d_n$ is an
explicit sum of {\it positive} integers. To the best of the authors'
knowledge those formulas have not appeared elsewhere.

\proclaim Proposition 6.2. The following expressions hold:
$$\leqalignno{\noalign{\vskip-4pt}
\qquad D_n(q)&=\kern-4pt\sum_{2\le 2k \le
n-1}\kern-4pt {1-q^{2k}\over 1-q}\,
{(q^{2k+2};q)_{n-2k-1}\over (1-q)^{n-2k-1}}\,
q^{{2k\choose 2}}+q^{{n\choose 2}}\chi(n\ {\sl even}),
&(6.18)\cr
d_n&=\sum_{2\le 2k\le n-1}(2k)(2k+2)_{n-2k-1}+\chi(n\ {\sl
even}).&(6.19)\cr
\noalign{\vskip-7pt}}
$$

{\it Proof.}\quad
When $q=1$, then (6.18) is transformed into (6.19). As for
(6.18), an easy $q$-calculation shows that its right-hand
side satisfies (6.9) when $Z=0$ and $q$ replaced
by~$q^{1/2}$.\qed

\medskip
Now, let ${}^{\inv\kern-2pt}A_n(q,Y_0):={}^{\ell\kern-1.5pt}B_n(q,Y_0,0,0)$.
Then
$$\leqalignno{
{}^{\inv\kern-2pt}A_n(q,Y_0)&=
\sum_{\sigma\in {\goth
S}_n}q^{\inv\sigma}\,Y_0^{\pix\sigma}\quad
(\pix:=\pix^+).&(6.20)\cr
\noalign{\hbox{Formula (1.7) specializes into}}
\sum_{n\ge 0}
{u^n\over (q;q)_n}\,{}^{\inv\kern-2pt}A_n(q,Y_0)
&=\Bigl(1-{u\over 1-q}\Bigr)^{-1} {(u;q)_\infty\over
(uY_0;q)_\infty};
&(6.21)\cr
\noalign{\hbox{In view of (6.15) we conclude that}}
A_n(q,Y_0)&={}^{\inv\kern-2pt}A_n(q,Y_0).&(6.22)\cr}
$$

For each permutation $\sigma=\sigma(1)\cdots \sigma(n)$ let
the {\it ligne of route} of~$\sigma$ be defined by
$\Ligne\sigma:=\{i:\sigma(i)>\sigma(i+1)\}$ and the {\it
inverse ligne of route} by
$\Iligne\sigma:=\Ligne\sigma^{-1}$. Notice that
$\maj\sigma=\sum_i i\,\chi(i\in\Ligne\sigma)$; we also let
$\imaj\sigma:=\sum_i i\,\chi(i\in\Iligne\sigma)$.
Furthermore, let ${\bf i}:\sigma\mapsto \sigma^{-1}$.
If~$\Phi$ designates the {\it second fundamental
transformation} described in [Fo68], [FoSc78], it is known
that the bijection $\Psi:={\bf i}\,\Phi\,{\bf i}$ of~${\goth
S}_n$ onto itself has the following property:
$(\Ligne,\imaj)\,\sigma=(\Ligne,\inv)\,\Psi(\sigma)$. Hence,
$$\leqalignno{\noalign{\vskip-5pt}
(\pix,\imaj)\,\sigma&=(\pix,\inv)\,\Psi(\sigma)
&(6.23)\cr
\noalign{\hbox{and then $A_n(q,Y_0)$ has the other
interpretation:}}
A_n(q,Y_0)&=\sum_{\sigma\in {\goth S}_n}q^{\imaj\sigma}
Y_0^{\pix\sigma}.&(6.24)\cr}
$$

Finally, let $K_n(q):=\sum\limits_{\sigma\in
K_n}q^{\inv\sigma}$. Then, with $Y_0:=0$ in (6.22) we have:
$$
K_n(q)={}^{\inv\kern-2pt}A_n(q,0)=A_n(q,0)=D_n(q).\leqno(6.25)
$$
However, it can be shown directly that $K_n(q)$ is equal to
the right-hand side of (6.18), because the sum occurring
in (6.18) reflects the geometry of the desarrangements. The running
term is nothing but the generating polynomial for the
desarrangements of order~$n$ whose leftmost trough is at
position~$2k$ by the number of inversions ``inv."

The bijection $\Psi$ also sends $K_n$ onto itself, so that
$$
\sum_{\sigma\in K_n}q^{\inv\sigma}
=\sum_{\sigma\in K_n}q^{\imaj\sigma},\leqno(6.26)
$$
a result obtained in this way by D\'esarm\'enien
and Wachs [DeWa90, 93], who also proved 
that for every subset $E\subset [n-1]$ we
have
$$
\#\{\sigma\in D_n:\Ligne\sigma=E\}
=\#\{\sigma\in K_n:\Iligne\sigma=E\}.\leqno(6.27)
$$

\vskip .5cm
\centerline{\bf References}
\medskip
{\eightpoint

\article ABR01|Ron M. Adin, Francesco Brenti and Yuval Roichman|Descent
Numbers and Major Indices for the Hyperoctahedral
Group|Adv. in Appl. Math.|27|2001|210--224|

\divers ABR05|Ron M. Adin, Francesco Brenti and Yuval Roichman|%
Equi-distribution over Descent Classes of the Hyperoctahedral 
Group, to appear in {\it J. Comb. Theory, Ser. A.}, {\oldstyle 2005}|

\article AR01|Ron M. Adin, Yuval Roichman|The flag major index
and group actions on polynomial rings|Europ. J.
Combin.|22|2001|431--6|

\livre An76|George E. Andrews|The Theory of
Partitions|Addison-Wesley, Reading MA, {\oldstyle 1976} ({\it
Encyclopedia of Math.and its  Appl., {\bf 2}})|

\livre Bo68|N. Bourbaki|Groupes et alg\`ebres de Lie, chap.~4, 5,
6|Hermann, Paris, {\oldstyle 1968}|

\article Br94|Francesco Brenti|$q$-Eulerian Polynomials
Arising from Coxeter Groups|Europ. J. 
Combin.|15|1994|417--441|

\article ClHaZe97|Robert J. Clarke, Guo-Niu Han, Jiang Zeng|A
Combinatorial Interpretation of the Seidel Generation of
$q$-derangement Numbers|Annals of Combin.|1|1997|313--327|

\divers Ch06|Chak-On Chow|On derangement polynomials of
type~$B$, {\it S\'em. Lothar. Combin.}, B55b
{\oldstyle 2006}, 6~p|

\divers ChXu06|William Y. C. Chen, Deheng Xu|Labeled Partitions and the
$q$-Derangement Numbers, {\it arXiv:math.CO/0606481} v1 20 Jun
{\oldstyle 2006}, 6~p|

\livre Co70|Louis Comtet|Analyse Combinatoire, vol.~2|Presses
Universitaires de France, Paris, {\oldstyle1970}|

\divers De84|Jacques D\'esarm\'enien|Une autre
interpr\'etation du nombre des d\'erange\-ments, {\sl S\'em.
Lothar. Combin.}, B08b, {\oldstyle 1982}, 6 pp. (Publ. I.R.M.A.
Strasbourg, {\oldstyle1984}, 229/S-08, p.~11-16)|

\divers DeWa88|Jacques D\'esarm\'enien, Michelle L.
Wachs|Descentes des d\'erangements et mots circulaires, {\sl
S\'em. Lothar. Combin.}, B19a, {\oldstyle1988}, 9 pp. (Publ.
I.R.M.A. Strasbourg, {\oldstyle1988}, 361/S-19, p.~13-21)|

\article DeWa93|Jacques D\'esarm\'enien, Michelle L.
Wachs|Descent Classes of Permutations with a Given Number of
Fixed Points|J. Combin. Theory, Ser.~A|64|1993|311--328|

\article Fo68|Dominique Foata|On the Netto inversion 
number of a sequence|Proc. Amer. Math. Soc.|19|1968|236--240|

\divers FoHa05|Dominique Foata, Guo-Niu Han|Signed words
and permutations, II; the Euler-Mahonian Polynomials, {\sl 
Electronic J. Combin.}, vol.~{\bf 11(2)}, {\oldstyle 2005},
\#R22, 18~p. (The Stanley Festschrift)|

\article FoSc78|Dominique Foata, Marcel-Paul
Sch\"utzenberger|Major index and inversion number of
permutations|Math. Nachr.|83|1978|143--159| 

\livre GaRa90|George Gasper,
Mizan Rahman|Basic Hypergeometric Series|London,
Cambridge Univ. Press, {\oldstyle 1990}  ({\sl Encyclopedia of
Math. and Its Appl.}, {\bf 35})|

\article GeRe93|Ira M. Gessel, Christophe Reutenauer|Counting
Permutations with Given Cycle Structure and Descent Set|J.
Combin. Theory, Ser.~A|64|1993|189--215|

\livre Hu90|James E. Humphreys|Reflection Groups and Coxeter
Groups|Cambridge Univ. Press, Cambridge (Cambridge Studies in
Adv. Math., {\bf 29}), {\oldstyle 1990}|

\article Re93a|V. Reiner|Signed permutation statistics|Europ. J.
Combinatorics|14|1993|553--567|

\article Re93b|V. Reiner|Signed permutation statistics and cycle
type|Europ. J. Combinatorics|14|1993|569--579|

\article Re93c|V. Reiner|Upper binomial posets and signed
permutation statistics|Europ. J.
Combinatorics|14|1993|581--588|

\article Re95a|V. Reiner|Descents and one-dimensional characters for
classical Weyl groups|Discrete Math.|140|1995|129--140|

\divers Re95b|V. Reiner|The distribution of descents and length
in a Coxeter group, {\sl Electronic J. Combinatorics}, vol.~{\bf 2},
{\oldstyle 1995}, \# R25|

\divers ReRo05|Amitai Regev, Yuval Roichman|Statistics on Wreath Products and 
Generalized Binomial-Stirling Numbers, to appear in {\it Israel J. Math.}, 
{\oldstyle 2005}|

\livre Ri58|John Riordan|An Introduction to Combinatorial
Analysis|New York, John Wiley \& Sons, {\oldstyle 1958}|

\article Wa89|Michelle L. Wachs|On $q$-derangement
numbers|Proc. Amer. Math. Soc.|106|1989|273--278|

}

\bigskip
\hbox{\vtop{\halign{#\hfil\cr
Dominique Foata \cr
Institut Lothaire\cr
1, rue Murner\cr
F-67000 Strasbourg, France\cr
\noalign{\smallskip}
{\tt foata@math.u-strasbg.fr}\cr}}
\qquad
\vtop{\halign{#\hfil\cr
Guo-Niu Han\cr
I.R.M.A. UMR 7501\cr
Universit\'e Louis Pasteur et CNRS\cr
7, rue Ren\'e-Descartes\cr
F-67084 Strasbourg, France\cr
\noalign{\smallskip}
{\tt guoniu@math.u-strasbg.fr}\cr}}
}

\bye